\documentstyle{amsppt} \magnification=\magstep1 \hsize=6 
truein \hcorrection{.375in} \vsize=8.5 truein 
\parindent=20pt \baselineskip=14pt \TagsOnRight
\NoBlackBoxes\footline{\hss\tenrm\folio\hss}

\centerline{A Survey of}
\centerline{Conformally Invariant Measures}
\centerline{on $H^m(\Delta )$}

\vskip1truein 

\centerline{Doug Pickrell}
\centerline{Department of Mathematics}
\centerline{University of Arizona}
\centerline{Tucson, Arizona 85721}
\centerline{pickrell\@math.arizona.edu}

\vskip1truein 

\flushpar Abstract.  The universal covering of $PSU(1,1)$ 
acts naturally on the space of holomorphic differentials 
of order $m$ on the unit disk.  The main purpose of this 
paper is to survey, as broadly as I am able, some basic 
sources and examples of invariant measures for this 
action.  A problem for the future is to determine, or at 
least to organize in some useful way, all of the 
invariant measures.

\bigskip\bigskip

\centerline{\S 0. Introduction}

\bigskip

The group $G=PSU(1,1)=\{g=\pm\left(\matrix a&b\\
\bar {b}&\bar {a}\endmatrix \right):\vert a\vert^2-\vert b\vert^2
=1\}$ acts 
on the open unit disk $\Delta\subset \Bbb C$ by linear fractional 
transformations, 
$$g:z\to\frac {az+b}{\bar {b}z+\bar {a}}.\tag 0.1$$
This identifies $G$ with the group of all conformal 
automorphisms of $\Delta$, or equivalently with the group of 
all orientation-preserving isometries of $\Delta$, equipped with 
the Poincare metric, $ds=\frac {\vert dz\vert}{1-\vert z\vert^2}$.  

Let $\tilde {G}$ denote the universal covering of $G$.  For our 
purposes a useful model is 
$$\tilde {G}=\{\tilde {g}=(\left(\matrix a&b\\
\bar {b}&\bar {a}\endmatrix \right),A)\in SU(1,1)\times \Bbb C:e^
A=a\},\tag 0.2$$
where if $\tilde {g}_3=\tilde {g}_1\tilde {g}_2$, then 
$$A_3=A_1+A_2+log(1+\frac {b_1\bar {b}_2}{a_1a_2}),\tag 0.3$$
and the covering short exact sequence is 
$$0\to \Bbb Z\to\tilde {G}\to G\to 0,\tag 0.4$$
where $n\in \Bbb Z$ maps to $(\left(\matrix -1&0\\
0&-1\endmatrix \right)^n,i\pi n)$, and $\tilde {g}$ above maps 
to $g=\pm\left(\matrix a&b\\
\bar {b}&\bar {a}\endmatrix \right)\in G$.  

For each real number $m$, the action of $G$ on $\Delta$ lifts to an 
action of $\tilde {G}$ on $\kappa^m$, the `$m$th power' of the canonical 
bundle on $\Delta$.  This induces a natural action of $\tilde {G}$ on the 
space of holomorphic sections of $\kappa^m$, the holomorphic 
differentials of degree $m$.  We will denote such a 
differential, and the space of all such differentials, by 
$f(z)(dz)^m$ and $H^m=H^m(\Delta )$, respectively.  By slight abuse 
of notation, given $\tilde {g}\in\tilde {G}$ as in $(0.2)$, we will write the 
action of $\tilde {g}$ on $H^m$ as 
$$\tilde {g}:f(z)(dz)^m\to f(\frac {\bar {a}z-b}{-\bar {b}z+a})(-
\bar {b}z+a)^{-2m}(dz)^m.\tag 0.5$$
The precise interpretation of $(-\bar {b}z+a)^{-2m}$ as a 
holomorphic function of $z\in\Delta$ is given by 
$$(a-\bar {b}z)^{-2m}=e^{-2mA}(1-\frac {\bar {b}}az)^{-2m};\tag 0.6$$
this is a well-defined holomorphic function, because 
$\vert a\vert^2-\vert b\vert^2=1$ implies $\vert\frac {\bar {b}}a
\vert\le 1$, and hence $\vert\frac {\bar {b}}az\vert <1$ for $z\in
\Delta$.  
Thus, for nonintegral $m$, the action (0.5) does depend 
upon $\tilde {g}\in\tilde {G}$, not simply $g$.  

The set of all $\tilde {G}$-invariant probability measures on $H^
m(\Delta )$ 
is a convex set, and the extreme points are the ergodic 
measures.  The naive futuristic problem which we pose 
is the following:  determine all, or at least organize in 
some useful way the interesting, ergodic $\tilde {G}$-invariant 
probability measures on $H^m(\Delta )$.  Our modest goal in this 
paper is simply to survey some basic sources and 
examples of such measures.  There are at least two 
cases which are of special interest:  $m=1$, because of 
its relation to the theory of classical Toeplitz and Hankel 
operators (\S 3 and \S 4), and $m=2$, because of its relation 
to universal Teichmuller space and conformal maps (\S 5).  
The case $m=0$ is closely related to the case $m=1$, 
because of the equivariant sequence 
$$0\to \Bbb C\to H^0(\Delta )@>{\partial}>>H^1(\Delta )\to 0.\tag 0.7$$

Invariant measures on $H^m(\Delta )$ are also related to
invariant measures on the space $\Cal C(\Delta )$ of 
(possibly infinite) configurations on $\Delta$. This 
correspondence arises via the
natural map
$$Z:H^m\setminus \{0\}\to \Cal C(\Delta ):\theta\to Z(\theta ),$$
where $Z(\theta )$ denotes the zero set. This seems to be a 
complete mystery at present.

The space $H^m$ is infinite dimensional, and $G$ is three 
dimensional.  As a consequence the support of an 
invariant measure can range from an orbit of the group 
(which will resemble the sphere bundle of a finite area 
Riemann surface), to something infinite dimensional.  
Measures of the former type, transitive measures, 
possibly generate all invariant measures, in the sense 
that by applying operations, such as convolution, all such 
measures are limits of these combinations.  We are 
primarily interested in these limits, and insightful 
interpretations which we might be able to attach to 
their supports, or to their (heuristic) `densities', 
or to their Fourier transforms.  

The action of $\tilde {G}$ on $H^m$, for $m>0$, is well-known to be 
unitary and irreducible on a dense subspace; we will 
denote this Hilbert subspace by $H^m\cap L^2$ (although this 
notation is potentially misleading when $0<m<1/2$; see 
\S 1B).  In this unitary context, at least as a starting 
point, it is perhaps useful to think of an invariant 
measure with `Lebesgue type support' (a notion which I 
do not know how to formalize) as heuristically having 
the form 
$$\frac 1{\Cal Z}\rho (\theta )dV(\theta ),\tag 0.8$$
where $\rho$ (the `density') is an invariant function of 
$\theta\in H^m$, and $dV$ denotes a fictitious Riemannian volume 
element corresponding to the invariant unitary structure.  

The basic examples are the Gaussians, where
$$\rho (\theta )=exp(-\frac 1{2t}\vert\theta\vert^2_{H^m\cap L^2}
).\tag 0.9$$
These measures can be characterized in many ways, and
they are invariant and ergodic 
with respect to the full unitary group of the Hilbert 
space $H^m\cap L^2$, and exhaust all such possibilities.  A 
venerable result of Irving Segal, which in some ways 
remains mysterious, asserts that these Gaussian 
measures are also $\tilde {G}$-ergodic ([Se]).  This is a fairly 
simple consequence of the obvious fact that the finite 
dimensional unitary representations of $\tilde {G}$ are trivial.  
But it seems hard to understand what ergodicity implies, 
in terms of the dynamics of $\tilde {G}$ acting on a random 
$m$-differential.  

One of the main points of this paper is to show that 
there are many other natural invariant functions which 
lead to invariant measures, via the heuristic expression 
$(0.8)$.  The simplest nonGaussian example is the following.  
If $\theta\in H^1$, let 
$$x(z)=\int\theta =x_1z+x_2z^2+..\in H^0/\Bbb C,\tag 0.10$$
as in the correspondence $(0.7)$.  The Hankel operator 
corresponding to $x$, or $\theta$, is represented by the infinite 
Hankel matrix 
$$B(x)=\left(\matrix .\\
x_n&.&\\
.&x_n&.\\
x_2&.&&\\
x_1&x_2&.&&x_n&.\endmatrix \right).\tag 0.11$$
This operator is Hilbert-Schmidt precisely when 
$\theta\in H^1\cap L^2$, and the action of $G$ on $x$, or $\theta$, intertwines 
with a natural action of $G$ on $B(x)$ by unitary 
conjugation (this is explained in \S 3).  Consequently the 
determinant 
$$D(\theta )=det(1+B(x)B(x)^{*})\tag 0.12$$
is a well-defined invariant function of $\theta\in H^1\cap L^2$.  For 
each $-1<l$, the heuristic expression 
$$\frac 1{\Cal Z}\frac 1{det(1+B(x)B(x)^{*})^{3+l}}dm(\theta ),\tag 0.13$$
has a rigorous interpretation as a limit of the finite 
dimensional probability measures 
$$\frac 1{\Cal Z_n}\frac 1{det(1+B_nB_n^{*})^{1+p_n+l}}dm(x_1,..,
x_n),\tag 0.14$$
where $B_n=B(\sum^nx_jz^{j-1})$, and by a miracle, $\Cal Z_n$and $
p_n$ can 
be calculated exactly.  This limit defines an invariant 
probability measure on $H^1$.  

This construction fits into a general pattern, which can 
be crudely formulated in the following way.  Suppose 
that a group acts in an `isometric way' on a dense 
subspace $X_{\infty}$ (with respect to a Riemannian or symplectic 
structure) of an infinite dimensional space $X$.  Suppose 
that $X_{\infty}$ is well-approximated by finite dimensional 
spaces 
$$X_1\subset X_2\subset ..\subset X_n\subset ..\subset X_{\infty}
\subset X.\tag 0.15$$
Suppose that $D$ is a positive invariant function on $X_{\infty}$ 
(such as $(0.12)$) such that for each $n$ 
$$\int_{X_n}\frac 1{D^p}dV_{X_n}<\infty\tag 0.16$$
for sufficiently large exponents $p$, and let $p_n$ denote the 
critical exponent, the smallest exponent such that $(0.16)$ 
holds for $p_n<p$.  Then with appropriate hypotheses 
one might expect that for $L>0$ 
$$\lim_{n\to\infty}\frac 1{\Cal Z_n}\frac 1{D^{p_n+L}}dV_{X_n}\tag 0.17$$
exists and defines an invariant measure on $X$. Although 
many examples are known, unfortunately a general 
existence result is missing.

In the Hankel example above, $p_n=2-\frac 1n$, which leads to 
the heuristic expression $(0.13)$.  In turn the expression 
$(0.13)$ suggests a potential interpretation (which I am 
unable to formalize) for the Fourier transform of the 
measure: up to a multiple, it is a fundamental 
solution for the 
pseudodifferential operator (involving infinitely many 
variables) $det(1+B(\frac {\partial}{\partial y})B(\frac {\partial}{
\partial y})^{*})^{3+l}$, where $y$ is dual to $x$.  

The plan of the paper is the following.  

In \S 1 we present elementary examples of ergodic 
invariant probability measures.  If $m$ is rational and 
nonnegative, then there exist transitive measures.  In 
addition to further discussing Gaussian measures, we 
also present an elementary example illustrating the 
general pattern outlined above.  

In $\S 2$ we discuss basic operations.  Because $H^m$ is a 
linear space, the convolution of two invariant measures 
is another invariant measure, and one can always scale a 
measure.  The multiplication map 
$$H^m\otimes H^n\to H^{m+n}:f(z)(dz)^m\otimes g(z)(dz)^n\to fg(dz
)^{m+n}\tag 0.18$$
is also $\tilde {G}$-equivariant, and we can use pushforward to 
define the `product' of two invariant measures. For 
example if $\nu^{(m)}_1$ denotes the standard Gaussian 
corresponding to $H^m\cap L^2$, as in $(0.8)$ and $(0.9)$, then the 
product $\nu^{(m)}_1\otimes\nu^{(n)}_1$ is an invariant measure on $
H^{m+n}$, and 
its Fourier transform is given by 
$$(\nu^{(m)}_1\otimes\nu^{(n)}_1\hat {)}(F)=\frac 1{det(1+B(m,n;F
)B(m,n;F)^{*})},\tag 0.19$$
where $B(m,n;F)$, for $F\in H^{m+n}\cap L^2$, is a kind of 
generalized Hankel operator which is represented by the 
matrix $(0.11)$ when $m=n=1/2$.  For both convolution and 
`multiplication', the product of two ergodic measures is 
another ergodic measure. Also the product 
operation intertwines via the zero map $Z$ 
with a natural convolution 
operation for measures on configuration space.

In \S 3 we discuss the measures in $(0.13)$, and 
analogues involving the operators $B(m,n;F)$, and  
more general linear equivariant operator-valued functions.  
At this point I do not know whether these measures are 
ergodic, whether they could possibly be generated from 
simpler measures by the operations described above, 
whether they remember the parameters $(m,n)$, and so on.  

In \S 4 I briefly discuss how $(0.13)$ likely fits into a 
larger context.  Given a simply connected compact 
symmetric space $U/K$, it is known that there exists a 
$LU$-invariant probability measure on a distributional 
completion of the loop space $L(U/K)$ (see [Pi2]).  This 
measure is believed to be reparameterization invariant.  
By expressing a loop in terms of its Riemann-Hilbert 
factorization, and by making use of a global version of 
the map $(0.7)$, 
$$H^0(\Delta ,U^{\Bbb C})\to H^1(\Delta ,\frak u^{\Bbb C}):g\to g^{
-1}\partial g,\tag 0.20$$
this should give rise to a conformally invariant measure 
on $H^1(\Delta ,\frak u^{\Bbb C})$, where $\frak u^{\Bbb C}$ denotes the Lie algebra of the 
complexification of $U$.  Hence the distributions of the 
matrix coefficients should provide a large source of 
$PSU(1,1)$-invariant measures on $H^1(\Delta )$.  In the special 
case in which $U/K$ is a group, which we consider in \S 4, 
there is an expression for the $H^1(\Delta ,\frak u^{\Bbb C})$-distribution 
which displays its conformal invariance.  The measures 
in $(0.13)$ appear to be related to the matrix coefficient 
corresponding to the highest root, when $U/K$ is $SU(2)$. 
But this merely looks plausible, and in particular, I do 
not have any insight into how to generate formulas for 
other matrix coefficients.  

In an analogous way invariant measures associated to 
$Diff(S^1)$, which are conjectured to exist, should give 
rise to $PSU(1,1)$-invariant measures on $H^2(\Delta )$.  This is 
related to an infinitesimal map 
$$H^{-1}(\Delta )@>{dS}>>H^2(\Delta )\tag 0.21$$
and a corresponding global map 
$$H^0(\Delta ,\Bbb C)@>{S}>>H^2(\Delta )\tag 0.22$$
where $S$ is the Schwarzian derivative, which are 
analogous to $(0.7)$ and $(0.20)$, respectively.  This 
speculation is described briefly in \S 5. I do not know if 
there is any useful connection with the Schramm-Loewner
evolution process.  

In the last section, \S 6, to compensate somewhat for the 
scattered character of the results in this paper, I have 
attempted to summarize some of the main points, and to 
pose some questions which point in promising directions.  
 
\bigskip

\centerline{Preliminaries}

\bigskip

The space $H^m$ is a Frechet space with respect to the 
topology of uniform convergence on compact sets, hence 
has a natural standard Borel structure.  With respect to 
this topology, the action $(0.5)$ is continuous.  

Throughout the remainder of this paper, the phrase `$\nu$ 
is an invariant measure' will be understood as an 
abbreviation for $`\nu$ is a $\tilde {G}$-invariant probability measure 
on $H^m$, relative to its natural standard Borel structure, 
for some $m$', unless explicitly stated otherwise.  

Given a Riemann surface $\Sigma$, $\kappa$, or for emphasis $\kappa
\vert_{\Sigma}$, will 
denote the canonical holomorphic line bundle.  A local 
holomorphic section of $\kappa$ is simply a holomorphic one 
form.  Thus $H^1(\Delta )=H^0(\kappa\vert_{\Delta})$, and in general 
$$H^m=H^m(\Delta )=H^0((\kappa\vert_{\Delta})^m).\tag 0.23$$
The line bundle $\kappa\vert_{\Delta}$ has a natural, i.e.  $\tilde {
G}$-invariant, 
Hermitian norm, where 
$$\vert (dz)^m\vert =(1-\vert z\vert^2)^m.\tag 0.24$$
This can be used to define various $\tilde {G}$-invariant Banach 
subspaces, such as $H^m\cap L^p$.  The global Poincare 
metric on $\Delta$ is given by 
$$d(z,w)=arctanh(\frac {z-w}{1-\bar {z}w}).\tag 0.25$$

\smallskip

\flushpar Configuration Space.

\smallskip

Let $\Cal C(\Delta )$ denote the configuration space of $\Delta$. As a 
topological space, $\Cal C(\Delta )$ divides up into connected 
components 
$$\Cal C(\Delta )=\bigsqcup_{0\le n\le\infty}\Cal C(\Delta )^{(n)}
,\tag 0.26$$
where $\Cal C(\Delta )^{(n)}$ denotes the space of configurations with $
n$ 
points. If $n<\infty$, then there is a universal covering
$$Perm(n)\times \{(q_i)\in\Delta^n:q_i\ne q_j,i\ne j\}\to \Cal C(
\Delta )^{(n)}.\tag 0.27$$
Suppose that $n=\infty$. By definition $\Cal C(\Delta )^{(\infty 
)}$ consists of 
countable subsets which have finite intersection with 
bounded sets. The topology is the weak topology 
generated by functions of the following form: given 
a compactly supported continuous function $f:\Delta\to \Bbb R$,
$$\Cal C(\Delta )^{(\infty )}:\to \Bbb R:\gamma\to\sum_{q\in\gamma}
f(q)\tag 0.28$$
is continuous. It is known that this topology is standard 
(defined by a complete separable metric), and 
the fundamental group is isomorphic 
to $Perm(\Bbb N)$. However this space is not locally simply 
connected, so that it does not have a universal covering.

For any $m$, there is a continuous equivariant map
$$Z:H^m(\Delta )\setminus \{0\}\to \Cal C(\Delta ):\theta\to Z(\theta 
),\tag 0.29$$
where $Z(\theta )$ is the zero set of $\theta$. For an invariant 
measure $\nu$ on $H^m$, $Z_{*}\nu$ will be an invariant measure on 
the configuration space, and if $\nu$ is ergodic, then $Z_{*}\nu$ 
will be ergodic. Thus an ergodic invariant measure
has an invariant, $n$, where $Z_{*}\nu$ is support on $\Cal C(\Delta 
)^{(n)}$. 

The simplest and most important measures on the 
infinite configuration space 
are 
the Poisson measures, $p_{\lambda}$, where $\lambda >0$. Given a finite 
disjoint collection of compact sets $K_i$, and integers $n_i\ge 0$, 
$$p_{\lambda}\{\gamma\in\Gamma_{\Delta}:\forall i,card(\gamma\cap 
K_i)=n_i\}=\prod_i\frac {(\lambda m(K_i))^{n_i}}{n_i!}e^{-\lambda 
m(K_i)}\tag 0.30$$
where $m$ denotes the invariant measure (see [KMM]).

\bigskip

\centerline{\S 1. Elementary Examples of Ergodic Measures.}

\bigskip

\flushpar A.  Transitive Measures.  

\smallskip

Suppose that $\nu$ is a transitive measure.  This means that 
there is a differential $\theta =f(z)(dz)^m\in H^m(\Delta )$ such that 
the orbit $\tilde {G}_{*}\theta$ has finite volume, with respect to the 
essentially unique $\tilde {G}$-invariant volume element on this 
space, and $\nu$ is the normalized volume element on this 
orbit.  

There are two possibilities for the stability subgroup of 
$\tilde {G}$ at $\theta$.  The first possibility is that it is the entire 
group, in which case $m=0$ and $\theta =f$ is constant, or 
$m>0$ and $\theta =0$.  In these cases $\nu =\delta_{\theta}$.  

The second possibility is that it is a lattice $\Cal L$, i.e.  a 
discrete subgroup of $\tilde {G}$ with a quotient which has finite 
volume.  The structure of such a lattice is known:  the 
projection of $\Cal L$ to $G$ is a lattice $L$, and $\Cal L\cap \Bbb Z
=N\Bbb Z$, for 
some positive integer $N$, where $\Bbb Z$ is identified with the 
center of $\tilde {G}$ as in $(0.4)$ (see Corollary 4.7.3, page 278 of 
[T]).  Let $G^{(N)}$ denote the unique $N$-fold covering of $G$ 
(for example $G^{(2)}=SU(1,1)$, $G^{(4)}$ is the $2$-fold covering 
of $SU(1,1)$, the metaplectic group, and so on).  On the 
one hand $\Cal L$ is the full inverse image of a lattice $L^{(N)}$ in 
$G^{(N)}$:  
$$\matrix 0&\to&N\Bbb Z&\to&\tilde {G}&\to&G^{(N)}&\to&0\\
&&\parallel&&\uparrow&&\uparrow\\
0&\to&N\Bbb Z&\to&\Cal L&\to&L^{(N)}&\to&0\endmatrix ;\tag 1.1$$
and on the other hand $L^{(N)}$ is isomorphic to $L$ with 
respect to the projection $G^{(N)}\to G$, so that it represents 
a splitting of the induced extension:  
$$\matrix 0&\to&\Bbb Z_N&\to&G^{(N)}&\to&G&\to&0\\
&&\uparrow&&\uparrow&&\uparrow\\
&&0&\to&L^{(N)}&\to&L&\to&0\endmatrix \tag 1.2$$

\flushpar(1.3) Examples:  (a) If $m$ is an integer, so that $G$ 
acts on $H^m$, then $\Cal L$ is simply the full inverse image of 
the lattice $L\subset G$.  As an example, $L$ could be the 
fundamental group of a compact Riemann surface $X$, 
where $X=L\backslash\Delta$.  In this context $\theta$ could be the pullback 
of a holomorphic one form, or quadratic differential, etc, 
from $X$ to $\Delta$.  

(b) Using a conformal isomorphism of $\Delta$ with the upper 
half plane, we can identify $G$ with $PSL(2,\Bbb R)$.  Let 
$L=\Gamma (N)$, the principal congruence subgroup of level $N$ 
(see page 54 of [FK]).  As a group $\Gamma (N)$ is free.  It 
follows that the extension of $L$ induced by 
$$\matrix 0&\to&\Bbb Z_N&\to&G^{(N)}&\to&G&\to&0\endmatrix \tag 1.4$$
can be split, and there is an essentially unique way to 
do the splitting.  This determines $L^{(N)}$ as in $(1.2)$, and $
\Cal L$ 
is the inverse image of $L^{(N)}$ in $\tilde {G}$.  In this context 
$$\theta =\eta (\tau )(d\tau )^{1/2},\text{\rm \ where }\eta (\tau 
)=q^{1/24}\prod_1^{\infty}(1-q^n)\tag 1.5$$
and $q=exp(2\pi i\tau )$, $0<Im(\tau )$, is a nonvanishing element of 
$H^{1/2}$ which is invariant with respect to $[\Gamma (1),\Gamma 
(1)]$, the 
commutator subgroup, which has index $6$.  By taking 
roots of $\eta$, it follows that there exist
transitive measures 
for arbitrary rational $m$, where the stability subgroup is 
commensurate with the $\Gamma (N)$.  

\smallskip

Returning to the question of stability subgroups, the 
other a priori possibilities are that the quotient is a 
circle or a disk.  But there are not invariant (finite) 
measures in these cases.  

We can summarize this discussion in the following way.  

\proclaim{(1.6) Proposition} Suppose $m>0$ and $\nu$ is a 
transitive $\tilde {G}$-invariant measure on $H^m(\Delta )$.  Then either 
$\nu =\delta_0$, or $m=\frac MN$ is rational (in reduced form), and $
\nu$ is 
the normalized invariant probability measure on the orbit 
$G^{(N)}_{*}\theta$, where $\theta$ is the pullback to $\Delta$ of a holomorphic 
$m$-differential on the finite volume Riemann surface $L\backslash
\Delta$, 
where $L\subset G$ is the image of $L^{(N)}$:  
$$\matrix G^{(N)}\times\kappa^{M/N}&\to&\kappa^{M/N}\\
\downarrow&&\downarrow\\
G\times\Delta&&L\backslash\Delta\endmatrix ;$$
here $\kappa^{1/N}$ denotes a fixed $N$th root of the canonical 
bundle $\kappa$ over $L\backslash\Delta$.  
\endproclaim

Given an invariant measure, it is natural to ask for 
minimal information which is needed to identify the 
measure.  At least for transitive measures this should 
be answered by the following 

\proclaim{ (1.7) Conjecture}The evaluation map 
$$eval_0:H^m(\Delta )\to \Bbb C:f(z)(dz)^m\to f(0)$$
induces an injective map from the set of transitive 
invariant probability measures to the set of 
$Rot$-invariant probability measures on $\Bbb C$.  
\endproclaim

In the case $m=0$, this asserts that a holomorphic 
function on $\Delta$, which is the pullback of a function on a 
finite volume surface, is determined by the (finite) 
distribution of its values, up to an automorphism of $\Delta$.  
This is analogous to the assertion that a function on a 
finite set is determined by a listing of its values and 
their multiplicities, up to a permutation of the 
function's domain.  

\bigskip

\flushpar B.  Gaussian Measures.  

\smallskip

For each $m\ge 0$, the action $\tilde {G}\times H^m(\Delta )$ contains an 
irreducible invariant unitary action.  If $m>\frac 12$, then the 
essentially unique invariant Hilbert space norm is given 
by (using the invariant norm $(0.24)$, and the invariant 
hyperbolic volume element) 
$$\vert f(dz)^m\vert^2_{H^m\cap L^2}=\frac {2m-1}{2\pi}\int_{\Delta}
\vert f(z)\vert^2(1-\vert z\vert^2)^{2m-2}dx\wedge dy=$$
$$\frac {2m-1}{2\pi}\sum_{n\ge 0}\vert f_n\vert^22\pi\int_0^1r^{2
n}(1-r^2)^{m-2}rdr=\sum_{n\ge 0}\vert f_n\vert^2\frac {\Gamma (n+
1)\Gamma (2m)}{\Gamma (2m+n)}$$
$$=\sum_{n\ge 0}\vert f_n\vert^2B(n+1,2m)=\vert f_0\vert^2+\frac 
1{2m}\vert f_1\vert^2+\frac {2!}{(2m)(2m+1)}\vert f_2\vert^2+..\tag 1.8$$
where $f=\sum f_nz^n$.  The seemingly unnatural factor $\frac {2m
-1}{2\pi}$ 
is inserted because the last sum shows that this norm 
can be analytically continued to $m>0$, since the Beta 
function coefficients satisfy $B(n+1,2m)>0$.  In the 
critical case $m=\frac 12$, the norm can be understood as the 
$L^2$ norm of the $S^1$-boundary values of $f$ with respect to 
$\vert dz\vert =d\theta$, which is manifestly $\tilde {G}$-invariant.  For 
$0<m<\frac 12$, there does not appear to be an integral 
representation for the norm, but we will nonetheless 
continue to denote the norm by $\vert\cdot\vert_{H^m\cap L^2}$.  This norm 
continues to be $\tilde {G}$-invariant (it would obviously be 
desirable, if it is possible, to find a geometric realization 
which clearly displays the invariance).  

The limit $m\downarrow 0$ is exceptional.  The action of $\tilde {
G}$ on $H^0$ is 
reducible:  $\Bbb C$ (the constants) is an invariant subspace, 
with quotient determined by $\partial :H^0(\Delta )\to H^1(\Delta 
)$.  The norm 
as defined by $(1.8)$ is not well-defined when $m=0$.  
However, in the definition of the norm, we can multiply 
by $m$.  In this case we obtain a Hilbert space 
substructure for the quotient $H^0/\Bbb C$ (the original vacuum, 
the constant $1$, is now a `ghost'), and an isometry 
$$(H^0/\Bbb C)\cap L^2@>{\partial}>>H^1\cap L^2:x\to\theta =\partial 
x.\tag 1.9$$
We will make extensive use of this natural isomorphism 
of $H^0/\Bbb C$ and $H^1$.  

To each of these unitary representations there is a 
corresponding one parameter convolution semigroup of 
invariant Gaussian measures, written heuristically as 
$$d\nu_T^{(m)}=\frac 1{\Cal Z_T}exp(-\frac 1{2T}\vert\theta\vert^
2_{H^m\cap L^2})dV(\theta ).\tag 1.10$$
These measures are ergodic, by a general result of 
Irving Segal ([Se]).  

In our context the Gaussian measure $(1.10)$ is realized as 
an infinite product measure 
$$\prod_{n\ge 0}\frac {B(n+1,2m)}{2\pi T}e^{-\frac 1{2T}B(n+1,2m)
\vert f_n\vert^2}dm(f_n).\tag 1.11$$
In other words we are given a sequence of independent 
standard Gaussian complex random variables $Z_n$, and we 
are considering the independent random variables 
$$f_n=a_nZ_n,\quad\text{\rm where }a_n=(\frac 1TB(n+1,2m))^{-1/2}
=(\frac {(2m)_nT}{n!})^{1/2},\tag 1.12$$
and a random Taylor series 
$$f(z)=\sum a_nZ_nz^n.\tag 1.13$$

If $T=1$, note that when $m=1/2$, then $f=\sum Z_nz^n$.  When 
$m=1$, then $\theta =fdz,$ where $f=\sum n^{1/2}Z_nz^n$, and 
$$x=\int_0^z\theta =\sum n^{-1/2}Z_nz^n.\tag 1.14$$
The $L^2$ norm of this on the circle diverges 
logarithmically, so this almost, but does not quite, have 
boundary values.  

The series $(1.13)$ is considered in great detail in [K], 
and here we merely note some of the main points.  

\proclaim{ (1.15) Proposition} Let $\nu =\nu_T^{(m)}$.  

(a) The radius of convergence of $f$ is 
$$(\limsup a_n^{1/n})^{-1}=1\quad a.e.\quad f\quad [\nu ],$$
i.e.  $\nu$ is supported on $H^m(\Delta )$.  

(b) For $a.e.$ $z\in S^1$, $(1.13)$ diverges, and $f$ assumes all 
values in $\Bbb C$, for $a.e.$ $f$ $[\nu ]$.  

Suppose that $m=1$, and consider $x(z)$ in $(1.14)$.  

(c) Fix an angle $\alpha$.  For all $a\in \Bbb C$, 
$$\liminf_{r\uparrow 1}\vert x(re^{i\alpha})-a\vert =0,\quad a.e.
\quad x\quad [\nu ],$$
and in particular the $x$ image of each ray is dense in $\Bbb C$, 
$a.e.$ $x$ $[\nu ]$.  

(d) Again fix the angle $\alpha$.  Then 
$$x(re^{i\alpha})=O(\sqrt {\rho (r)ln(\frac {\rho (\sqrt r)}{1-r}
)})\quad as\quad r\uparrow 1,$$
$a.e.$ $x$ $[\nu ]$, where $\rho (r)=-ln(1-r^2)$.  

\endproclaim

For (a), see section 2, for $(b)$, see Theorem 1 of chapter 
13, for (c), see Theorem 3, page 184, and for (d), see $(22)$ 
on page 187 of [K].  

The conditions $(c)$ and $(d)$ give us some feeling for what 
a random $x(\theta )$ is like for the $H^1$-Gaussian.  It should be 
possible to refine these conditions to give criteria for a 
given $x$ to have a dense $G$-orbit (see $(6.3)$ below).  

The central limit theorem implies the following 

\proclaim{ (1.16) Proposition}Suppose that $\nu$ is an 
invariant measure with the property that $\phi_{*}\nu$ has finite 
variance for each $\phi\in (H^m\cap L^2)^{*}$.  Then 

(a) The variance of $\phi_{*}\nu$ is independent of $\phi$, provided 
$\vert\phi\vert_{H^m\cap L^2}=1$; denote this common value by $T$; and 

(b) the scaled $n$-fold convolution 
$$\nu *\nu *..*\nu (\sqrt n\cdot )\to d\nu_T^{H^m}\quad\text{\rm as }
\quad n\to\infty$$
in a weak sense.  
\endproclaim

\demo{Proof of (1.16) } Let $\nu_n$ denote the scaled $n$-fold 
convolution:  $\nu_n(E)=\nu^{n*}(\sqrt nE)$.  The invariance of $
\nu$ 
implies that the mean of $\nu$ is zero.  Given $\phi\in (H^m\cap 
L^2)^{*}$, 
let $\sigma^2=\sigma^2(\phi_{*}\nu )$, the variance of the $\phi$ distribution of $
\nu$.  
This is finite, by assumption.  The central limit theorem 
implies that 
$$\phi_{*}\nu_n\to\frac 1{\sigma\pi}e^{-\frac 1{2\sigma^2}\vert\phi
\vert^2}dm(\phi )\tag 1.17$$
Given $\phi$ and $\psi$, this implies that 
$$\int\phi\bar{\psi }d\nu_n\to\langle\phi ,\psi\rangle\quad\text{\rm as}
\quad n\to\infty ,\tag 1.18$$
where $\langle\cdot ,\cdot\rangle$ defines an invariant Hermitian inner 
product.  By irreducibility, this inner product is 
necessarily a multiple of $\vert\cdot\vert_{H^m\cap L^2}$.  This determines $
T$ 
and proves $(a)$.  

Given $\phi$ and $\psi$ of unit length with $\phi\perp\psi$, by considering 
the distribution for $s\phi +t\psi$, we see that in the limit 
$n\to\infty$, $\phi$ and $\psi$ are independent.  This implies that the 
pair $\phi ,\psi$ has a Gaussian distribution.  This implies $(b)$.  
\qed
\enddemo

Note that if $\nu$ is a transitive measure, then $\nu$, and its 
scalings, have finite variance, and hence the Gaussians 
are in the closure of the convolution algebra generated 
by $\nu$ and its scalings.  

For later use we recall the $n$-point functions for the 
standard Gaussians.  Suppose that we fix a point $z_0\in\Delta$.  
The continuous linear functional 
$$eval_{z_0}:H^m(\Delta )\to \Bbb C:f(z)(dz)^m\to f(z_0),\tag 1.18$$
restricted to `$L^2$' differentials, is represented by the 
differential $(1-\bar {z}_0z)^{-2m}(dz)^m$, that is, 
$$f(z_0)=\langle f(dz)^m,\frac {(dz)^m}{(1-\bar {z}_0z)^{2m}}\rangle_{
H^m\cap L^2}.\tag 1.19$$
Given $n$-points $z_i\in\Delta$, the covariance matrix $C$ has 
entries given by 
$$C_{ij}=\langle eval_{z_i},eval_{z_j}\rangle =(1-\bar {z}_iz_j)^{
-2m}.\tag 1.20$$
The $\nu_1^{(m)}$-distribution for the $n$-point function 
$$H^m(\Delta )@>{eval_{(z_1,..,z_n)}}>>\Bbb C^n:f(z)(dz)^m\to (f(
z_i))_{i=1,..,n}\tag 1.21$$
is given by the Gaussian measure 
$$2^{-n}det_{\Bbb C}C^{-1}exp(-\frac 14(C^{-1})_{jk}z_j\bar {z}_k
)\prod_1^ndm(z_j).\tag 1.22$$

Suppose that $m=1$.  In this case it is more to the point 
to consider the distribution 
$$H^1\to \Bbb C^n:\theta =f(z)dz\to (x(z_i)-x(0))_{1\le i\le n},\tag 1.23$$
where $x\in H^0/\Bbb C$ is related to $\theta$ as in $(1.14)$.  Now 
$$x(z_0)-x(0)=\sum x_nz_0^n=\sum nx_n(\frac {\bar {z}_0^n}n)^{*}=
\langle\sum x_nz^n,\sum (\frac {\bar {z}_0^n}n)z^n\rangle_{H^0/\Bbb C
\cap L^2}$$
$$=\langle x,ln(\frac 1{1-\bar {z}_0z})\rangle\tag 1.24$$
Therefore 
$$C_{ij}=\langle ln(\frac 1{1-\bar {z}_iz}),ln(\frac 1{1-\bar {z}_
jz})\rangle =ln(\frac 1{1-\bar {z}_iz_j}).\tag 1.25$$

\bigskip

\flushpar C.  Another Example.  

\smallskip

The following elementary example is included to 
illustrate what is involved in constructing and analyzing 
an invariant measure, and the general scheme outlined in 
the introduction.  

For simplicity suppose $m=1$.  Write $\theta\in H^1$ as 
$$\theta =(\theta_1+\theta_2z+..)dz,\tag 1.26$$
and define the projection $P_n\theta =(\theta_1+..+\theta_nz^{n-1}
)dz$.  Define 
a probability measure on $P_nH^1$ by 
$$d\mu_l^{(n)}=\frac 1{\Cal Z}(1+\vert P_n\theta\vert_{H^1\cap L^
2})^{-n-1-l}dm(P_n\theta )=$$
$$\frac {\Gamma (n+l+1)}{\pi^n\Gamma (l+1)}(1+\vert\theta_1\vert^
2+\frac 12\vert\theta_2\vert^2+..+\frac 1n\vert\theta_n\vert^2)^{
-n-1-l}dm(P_n\theta ).\tag 1.27$$
This measure is finite precisely when $-1<l$, for each $n$.  
In the context of the general scheme discussed in the 
introduction, 
$$D=1+\vert\theta\vert^2_{H^1\cap L^2}=1+tr(B(x)B(x)^{*}),\tag 1.28$$
and the critical exponent is $p_n=n$.  
  
These measures are coherent in the sense that 
$$(P_n)_{*}\mu_l^{(n)}=\mu_l^{(n-1)}.\tag 1.29$$
The formal completion of $H^1$ is the infinite product 
space which we write formally as 
$$H^1_{formal}=\{\theta =(\sum_1^{\infty}\theta_jz^{j-1})dz:(\theta_
j)\in\prod_1^{\infty}\Bbb C\}.\tag 1.30$$
The coherence property, together with Kolmogorov's 
theorem (page 228 of [B]), implies that there 
is a unique probability measure $\mu_l$ on the formal 
completion such that 
$$(P_n)_{*}\mu_l=\mu_l^{(n)},\tag 1.31$$
for all $n$.  There are a number of slightly different 
ways of showing that $\mu_l$ is supported on $H^1$:  (1) it is a 
convex combination of Gaussians ($(1.34)$ below), (2) one 
can compute its Fourier transform in spherical 
coordinates and observe that it is continuous on the dual 
of $H^1$ (see $(1.35)$), and (3) in the next section we will see 
that $\mu_l$ is a `quotient of Gaussians'.  

Along the lines of $(0.8)$, the measure $\mu_l$ has a heuristic 
expression 
$$d\mu_l=\frac 1{\Cal Z}\frac 1{(1+\vert\theta\vert^2_{H^1\cap L^
2})^{1+p+l}}dm(\theta ),\tag 1.32$$
where $p=\infty$ in this case.  Although of limited use, this 
expression does suggest that $\mu_l$ depends only upon the 
unitary structure of $H^1$, and hence that $\mu_l$ should be a 
convex combination of Gaussians.  In fact, for each $n$ 
$$\frac 1{\pi^n}\frac {\Gamma (n+1+l)}{\Gamma (1+l)}\frac 1{(1+\vert 
P_n\theta\vert^2)^{1+n+l}}dm(P_n\theta )=$$
$$\int_0^{\infty}\{(\frac {\beta}{\pi})^ne^{-\beta\vert P_n\theta
\vert^2}dm(P_n\theta )\}\frac 1{\Gamma (1+l)}\beta^le^{-\beta}d\beta\tag 1.33$$
and hence 
$$d\mu_l=\int_{\beta =0}^{\infty}\left(d\nu_{\beta^{-1}}^{(1)}\right
)\frac 1{\Gamma (1+l)}\beta^le^{-\beta}d\beta .\tag 1.34$$
These expressions are very interesting in the way that 
normalization constants and exponents interact.  

Note that this implies that $\mu_l$ is not ergodic, a fact 
which is not so easy to understand in a naive way.  

In general, given an invariant measure $\nu$ on $H^m$, for $m>0$, we 
define its Fourier transform using the unitary 
substructure, 
$$\hat{\nu }(F)=\int e^{-i\Bbb Re\langle F,f\rangle_{H^m\cap L^2}}
d\nu (f)\tag 1.35$$
for $F\in H^m_{alg}$ (this is well-defined for any finite measure 
on $H^m_{formal}$, because $H^m_{alg}$ is dual to the formal 
completion).  

It is interesting to note that the Fourier transform of 
$\mu_l$ can be written in two ways 
$$\int e^{-i\Bbb Re\langle\theta ,\theta'\rangle_{H^1}}d\mu_l(\theta 
)=\int_0^{\infty}\frac 1{(1+r^2)^2}J_0(r\vert\theta'\vert_{H^1\cap 
L^2})rdr$$
$$=\int_0^{\infty}e^{-\vert\theta'\vert^2/T}\frac 1{\Gamma (1+l)}
T^le^{-T}dT.\tag 1.36$$
The first expression is obtained by integrating with 
respect to polar coordinates; the second expression 
follows from $(1.34)$.  If we think of this Fourier 
transform as a function of $n$ variables, then it is a 
multiple of the fundamental solution of $(1+\Delta )^{1+n+l}$.  

Note that both expressions clearly show that the Fourier 
transform is defined and continuous on $(H^1)^{*}$.  Hence the 
support properties of $\mu_l$ and the Gaussians are roughly
the same.  

\bigskip

\centerline{\S 2. Operations and Further Examples.}

\smallskip

There are a number of basic operations which we can 
perform to obtain new invariant measures from known 
invariant measures.  

Because $\tilde {G}$ acts linearly on $H^m$, we can scale a given 
measure, we can form the convolution of two invariant 
measures, and we can form convex combinations.  

The multiplication map 
$$M:H^m\times H^n\to H^{m+n}:(\theta ,\theta')\to\theta\theta'=f(
z)g(z)(dz)^{m+n}\tag 2.1$$
is also $\tilde {G}$-equivariant, where $\theta =f(z)(dz)^m$ and 
$\theta'=g(z)(dz)^m$.  Hence if $\nu$ and $\mu$ are invariant measures 
on $H^m$ and $H^n$, respectively, then $M_{*}(\nu\times\mu )$ is an 
invariant measure on $H^{m+n}$.  We will refer to this 
measure as the multiplicative image of $\nu$ and $\mu$, and we 
will write $\nu\otimes\mu =M_{*}(\nu\times\mu )$.  

We can also consider the quotient map 
$$Q:(\theta ,\theta')\to\theta /\theta'=\frac {f(z)}{g(z)}(dz)^{m
-n},\tag 2.2$$
which is defined, provided $\theta'$ is nonvanishing.  If $\mu$ is 
supported on nonvanishing differentials, then $Q_{*}(\nu\times\mu 
)$ 
will be a well-defined invariant measure on $H^{m-n}$.  

It is definitely the case that there exist transitive $\mu$ on 
$H^n$ which are supported on nonvanishing differentials, 
assuming that $n$ is rational; this was one of main points 
of example $(b)$ of $(1.3)$.  However for measures having 
infinite dimensional support, one expects that the sample 
properties will be similar to the Gaussian case.  For the 
Gaussians, provided $m>0$, part $(b)$ of $(1.15)$ implies that 
with probability one, a typical $f$ vanishes at some 
points.  

For this reason, in thinking about $Q$, it is useful to 
consider the formal completions of the spaces $H^m$.  The 
formal completion of $H^m$, denote $H^m_{formal}$, is the space 
consisting of differentials $\theta =f(dz)^m$, where $f$ is simple 
a formal power series (see $(1.30)$).  We view $Q$ as a map 
$$Q:H^m_{formal}\times \{\theta'\in H^n_{formal}:g_0\ne 0\}\to H^{
m-n}_{formal}\tag 2.3$$
where $f/g$ in $(2.2)$ is interpreted as the formal power 
series 
$$\frac fg=\frac {f_0}{g_0}+\frac {f_1g_0-f_0g_1}{g_0^2}z+..\tag 2.4$$

\proclaim{ (2.5) Proposition } If $\nu$ and $\mu$ are ergodic 
measures, then, whenever defined, the convolution $\nu *\mu$, 
the multiplicative image $\nu\otimes\mu$, and the quotient image 
$Q_{*}(\nu\times\mu )$ are ergodic invariant measures.  
\endproclaim

This follows from the fact that $\tilde {G}$ does not have any 
finite dimensional unitary representations, see [Se].  

Now recall the natural map $H^m\to \Cal C(\Delta )$. There is a natural 
operation of convolution for measures on $\Cal C(\Delta )$, 
corresponding to the natural union operation
$$\Cal C(\Delta )\times \Cal C(\Delta )\to \Cal C(\Delta ):\gamma_
1,\gamma_2\to\gamma_1\cup\gamma_2.$$
It is well-known, for example, that for the Poisson 
measures
$$p_{\lambda_1}*p_{\lambda_2}=p_{\lambda_1}*p_{\lambda_2}.$$

\proclaim{ (2.6) Proposition} If $\lambda >0$, then $Z_{*}(\lambda_{
*}\nu )=Z_{*}(\nu )$.
Also $Z_{*}$ intertwines $\otimes$ and convolution.
\endproclaim

In this paper we will not make any use of this 
proposition, simply because we do not know how to 
compute $Z$ for any nontransitive examples.

\bigskip

\flushpar Some Calculations.  

\smallskip

We would like to be able to calculate the Fourier 
transform of the measure 
$$\nu^{(m_1)}_1\otimes ..\otimes\nu^{(m_n)}_1\tag 2.7$$
We are able to do this in a completely satisfactory way 
only for pairs.  We will need the fact that there is a 
essentially unique $\tilde {G}$-equivariant embedding 
$$H^{m+n}\cap L^2\to (H^m\cap L^2)\otimes (H^n\cap L^2)=$$
$$\Cal L_2(H^m(\Delta^{*})\cap L^2,H^n\cap L^2),\tag 2.8$$
and this embedding is completely determined by requiring 
that it map $(dz)^{m+n}$ to $(dz)^m\otimes (dz)^n$.  

The is a consequence of the following standard result.

\proclaim{(2.9)Lemma}As a $\tilde {G}$-representation 
$$(H^m\cap L^2)\otimes (H^n\cap L^2)=(H^{m+n}\cap L^2)\oplus (H^{
m+n+1}\cap L^2)\oplus ..\tag 2.10$$
Furthermore multiplication, $M$ in $(2.1)$, 
is the orthogonal projection 
$$0\to Ker\to (H^m\cap L^2)\otimes (H^n\cap L^2)@>{M}>>(H^{m+n}\cap 
L^2)\to 0.\tag 2.11$$

\endproclaim

\demo{Proof of (2.9)} The action of (the covering in $\tilde {G}$) of 
rotations extends to a holomorphic contraction 
representation on $H^m\cap L^2$.
The corresponding character of $H^m$ is $q^m/(1-q)$, where 
$q=e^{2\pi i\tau}$, $Im(\tau )>0$, and this character is a complete 
invariant of this (lowest weight) representation.  
Now simply observe that the characters of the two 
sides of $(2.10)$ are equal:  
$$\frac {q^m}{1-q}\cdot\frac {q^n}{1-q}=\frac {q^{m+n}+q^{m+n+1}+
..}{1-q}.\tag 2.12$$
The second claim follows from irreducibility of $H^{m+n}$.  
\qed
\enddemo

Given $F\in H^{m+n}\cap L^2$, we will denote the Hilbert- Schmidt 
operator corresponding to $F$ via $(2.8)$ by $B(m,n;F)$.  We 
will write out the matrix for this operator in $(2.17)$ 
below.  This calculation will make it clear that $B(m,n;F)$ 
is a kind of generalized Hankel operator (note that we 
are also using $B$ for the Beta function, but the number 
of arguments and context should always make it clear 
which is intended).  

\proclaim{(2.13) Proposition} 
$$(\nu^{(m)}_S\otimes\nu^{(n)}_T)\hat {(}F)=\frac 1{det(1+STB(m,n
;F)B(m,n;F)^{*})}$$
for $F\in H^{m+n}\cap L^2$.  
\endproclaim

\demo{Proof of $(2.13)$}.  Suppose $f=f(z)(dz)^m\in H^m$, 
$g=g(z)(dz)^n\in H^n$, and $F=F(z)(dz)^{m+n}\in H^{m+n}_{alg}$.  Then 
$$\langle F,fg\rangle_{H^{m+n}}=\sum_l\sum_{k\le l}F_l\bar {f}_{l
-k}\bar {g}_kB(l+1,2(m+n))$$
$$=\sum_k\left(\sum_jF_{k+j}\bar f_j\frac {B(k+j+1,2m+2n)}{B(k+1,
2n)}\right)\bar {g}_kB(k+1,2n).\tag 2.14$$
This last sum is the $H^n$ inner product of $f$ and $h$, where 
$h_k$ is the term in braces in $(2.14)$ (which depends upon 
both $F$ and $f$).  The equivariant properties of the linear 
mapping, depending upon $F$, which sends $\bar {f}$ to $h$ are the 
same as those for the operator $B(m,n;F)$.  The 
uniqueness of the embedding $(2.8)$ implies that 
$$\langle F,fg\rangle_{H^{m+n}}=\langle B(m,n;F)\bar {f},g\rangle_{
H^n}.\tag 2.15$$

\flushpar\ (2.16) Remark.  The calculation $(2.14)$ implies 
that the matrix of $B(m,n;F)$, relative to the orthogonal 
(but not necessarily orthonormal) bases $z^k(dz)^m$ 
and $z^l(dz)^n$ for $(H^m\cap L^2$$)^{*}$ and $H^n\cap L^2$, 
respectively, has entries given by 
$$B(m,n;F)_{j,k}=F_{k+j}\frac {B(k+j+1,2m+2n)}{B(k+1,2n)}\tag 2.17$$
If we scale these bases to obtain orthonormal bases, 
then the matrix entries are 
$$\frac {B(k+j+1,2m+2n)}{\sqrt {B(j+1,2m)B(k+1,2n)}}.\tag 2.18$$
In the special case that $m=n=\frac 12$, $B(\cdot ,2m)=B(\cdot ,2
n)=1$, 
$B(k+j+1,2m+2n)=\frac 1{k+j+1}$, and hence $B(\frac 12,\frac 12;F
)=B(x)$, 
where $x$ is the antiderivative of $F(z)dz$.  

\smallskip

Using $(2.15)$, we now see, by first integrating with 
respect to $g$, that 
$$(\nu_S^{(m)}\otimes\nu_T^{(n)})^{\hat{}}(F)=\int e^{-\frac T2\langle 
A\bar {f},f\rangle}d\nu_S^{(m)}(f),\tag 2.19$$
where $A=B(m,n;F)B(m,n;F)^{*}$.  There is a standard 
formula for the Gaussian integral of the exponential of a 
quadratic functional, as in $(2.19)$, and this implies $(2.13)$.  
\qed
\enddemo

\proclaim{(2.20) Corollary} (a) $det(1+B(m,n;F))B(m,n;F)^{*})^{-j}$ 
is a positive definite function of $F\in H^{m+n}\cap L^2$, for 
$j=1,2,..$.  

(b) The Fourier transform of $\nu_S^{(m)}\otimes\nu_T^{(n)}$ is continuous on 
$H^{m+n}\cap L^2$.  
\endproclaim

It is definitely interesting to ask whether (a) is true for 
nonintegral $j>0$ (and what kind of measure would 
correspond to this Fourier transform).  It is not so hard 
to see that (b) holds for any of the measures $(2.7)$.  But 
whether there is a tractable formula for the Fourier 
transform of $(2.7)$, in general, is unclear.  

Now we consider the quotient operation on formal 
completions.  

Here is a particular application of this construction 
which shows that in general an invariant probability on 
the formal completion of $H^m$ is not determined by its 
one-point function (which we conjectured is the case for 
transitive measures, in $(1.7)$).  

\proclaim{(2.21) Proposition}Fix $T>0$, and $m\ge 1$.  The 
image measure 
$$\nu_m=Q_{*}(\nu_T^{H^m}\times\nu_T^{H^{m-1}})$$
is a $G$-invariant probability measure on $H^1(\Delta )$ with the 
property that 
$$(eval_0)_{*}\nu_m=\frac 1{\Cal Z}\frac {dm(q_0)}{(1+\vert q_0\vert^
2)^2}.$$
\endproclaim

\demo{Proof of (2.21)} For all $m$, 
$(eval_0)_{*}\nu_T^{H^m}=\frac 1{\Cal Z}e^{-\frac 1{2T}\vert q_0\vert^
2}dm(q_0)$.  The proposition now 
follows from the commutativity of the diagram 
$$\matrix H^m\times H^{m-1}&@>{Q}>>&H^1\\
\downarrow eval_0&&\downarrow eval_0\\
\Bbb C\times \Bbb C&\to&\Bbb C&:(z,w)\to\frac zw\endmatrix \tag 2.22$$
together with the fact that a ratio of independent 
Gaussians (with the same temperature) is the measure 
above.  

\qed
\enddemo

\flushpar(2.23) Example.  The simplest example is when 
$m=1$, where we interpret $H^0$ to be $\Bbb C$ with the usual 
Gaussian.  In this case we obtain the measure on $H^1$ 
discussed in \S 1C.  This explains in another way why that 
measure is supported on $H^1$.  

\smallskip

Although a little tedious, it is a straightforward matter 
to compute the $n$-point function for $\nu_m$ in general (using 
the well-known fact that the $n$-point functions for a 
Gaussian are determined by the $2$-point functions (and 
note that when we say $n$-point function, we actually 
have $n$ complex values, so that we have $2n$ real values).  
It is given by 
$$\frac 1{2^ndetC(m)C(m-1)}\sum_{i,j\in Perm(n)}\prod_{l=1}^n[(1+
C(m-1)diag(q_k)C(m)^{-1}diag(\bar {q}_k))^{-1}C(m-1)]_{i_lj_l}$$
where $C(m)$ denotes the covariance for the $n$-point 
function corresponding to the points $q_k$ for the Gaussian 
$\nu_1^{H^m}$.  From this expression, say by looking at the case 
$n=2$, one can see that the $\nu$$_m$ are distinct measures, as 
$m$ varies.  Thus the one-point function certainly does 
not determine an invariant probability (on the formal 
completion).  

Using this construction it is possible to construct 
invariant probabilities on $H^m$ with many different 
one-point functions.  But given a candidate for a 
one-point function, i.e.  a $Rot$-invariant probability on $\Bbb C$, 
I have no idea how to decide whether it comes from an 
invariant probability on $H^m$.  

\bigskip

\centerline{\S 3. Linear Operator-Valued Functions, }
\centerline{and Associated Measures}

\bigskip

In the first part of this section, 
after some further discussion of the 
sense in which the map from $\theta\in H^1\cap L^2$ to the Hankel 
matrix $(0.11)$ is $PSU(1,1)$-equivariant, we will discuss the 
limit $(0.13)$.  In the second part of the section, using the 
operators $B(m,n;F)$ of $\S 2$ (and more general, generalized 
Hankel operators), we indicate one type of generalization 
(which also involve more 
complexity).    

\bigskip

\flushpar Hankel Operators and $m=1$.

\smallskip

The circle has a unique nontrivial (or Moebius) spin 
structure:  the ($C^r$) sections are of the form $f(\theta )(d\theta 
)^{1/2}$, 
where $f$ is a $\underline {real-valued}$ ($C^r$) function satisfying 
$f(\theta +2\pi )=-f(\theta )$ (Note that such an odd real-valued 
function must vanish, and this explains why the bundle 
is trivial).  The complexification of this real line is the 
restriction to $S^1$ of $\kappa^{1/2}$, the unique holomorphic square 
root of the canonical line bundle $\kappa$ for $\hat {\Bbb C}$.  This 
complexified bundle is trivial:  a global section is 
$$(dz)^{1/2}=\sqrt ie^{i\theta /2}(d\theta )^{1/2},\tag 3.1$$
where we fix a choice of square root for $i$.  The group 
$SU(1,1)$ (the double cover of $G$!)  acts on this bundle, as 
in $(0.5)$.  We will denote the space of sections of this 
complexification by $\Omega^{1/2}(S^1)$ (or $\Omega^0(\kappa^{1/2}
\vert_{S^1})$).  There is a 
natural Hermitian inner product 
$$\Omega^{1/2}(S^1)\otimes\Omega^{1/2}(S^1)^{conj}\to \Bbb C:\theta
\otimes\bar{\eta}\to\int_{S^1}\theta\bar{\eta}\tag 3.2$$
The point is that if $\theta =f(z)(dz)^{1/2}$ and $\eta =g(z)(dz)^{
1/2}$ 
are odd spinors on $S^1$, then $\theta\bar{\eta }=f\bar {g}\vert 
dz\vert$ is naturally a 
one density on $S^1$, which can be integrated.  We will 
write $\Omega^{1/2}(S^1)\cap L^2$ when we are thinking of this space 
as a Hilbert space.  

There is a Hilbert space isomorphism 
$$L^2(S^1,\Bbb C)\to\Omega^{1/2}(S^1)\cap L^2:f(z)\to f(z)(dz)^{1
/2}.\tag 3.3$$
In this identification the Hardy polarization of $L^2(S^1,\Bbb C)$ 
used in [PS] is identified with the $SU(1,1)$-equivariant 
polarization 
$$\Omega^{1/2}(S^1)\cap L^2=H^{1/2}\cap L^2\bigoplus H^{1/2}(\Delta^{
*})\cap L^2.\tag 3.4$$
The identification $(3.3)$ will be fixed from this point 
onward.  

As in [PS], given $x\in L^{\infty}(S^1,\Bbb C)$, we obtain a multiplication 
operator $M_x$ on $\Omega^{1/2}(S^1)$.  With respect to the Hardy 
polarization, we write 
$$M_x=\left(\matrix A&B\\
C&D\endmatrix \right),\tag 3.5$$
where $A$ is referred to as the Toeplitz operator 
associated to $g$, and $B=B(x)$ (or $C$) is the Hankel 
operator associated to $x$.  Relative to the basis 
$z^j(dz)^{1/2}$, the matrix of $M_x$ is constant along diagonals, 
and in particular the matrix of $B$ is given by $(0.11)$.  

Note that if $x\in H^0/\Bbb C$, then the $L^2$ norm of $x$ is the 
same as the Hilbert-Schmidt norm of $B(x)$.  These 
considerations imply the following 

\proclaim{(3.6)Proposition}The map $x\to B(x)$ is equivariant 
with respect the natural unitary action of $PSU(1,1)$ on 
$x\in H^0/\Bbb C\cap L^2$ and the action by unitary conjugation of 
$SU(1,1)$ on 
$$B(x)\in \Cal L_2(H^{1/2}(\Delta^{*})\cap L^2,H^{1/2}\cap L^2)=H^{
1/2}\cap L^2\otimes H^{1/2}\cap L^2.$$
\endproclaim

This explains why the determinant $(0.11)$ is 
$PSU(1,1)$-invariant, and why we expect the measure with 
formal expression $(0.13)$ to be $PSU(1,1)$-invariant.  

Now fix $N$.  Define $B_N(x)$ to be the finite Hankel matrix 
$$B_N(x)=\left(\matrix x_N&0&..\\
.&x_N&0&..\\
.&.&.\\
x_2&x_3&..\\
x_1&x_2&x_3&.&x_N\endmatrix \right)\tag 3.7$$

It is clear that for sufficiently large $p$, 
$$\frac 1{det(1+B_NB_N^{*})^p}\in L^1(dm(x_1,..,x_N)).\tag 3.8$$
Let $p_N$ denote the infimum of these $p$, and 
$$\Cal Z(p,N)=\int\frac 1{det(1+B_NB_N^{*})^p}dm(x_1,..,x_N)<\infty 
,\quad p_N<p.\tag 3.9$$
It is clear that $p_1\le p_2\le ..$.  

\smallskip

\flushpar(3.10) Examples.  (1) For $N=1$ 
$$\Cal Z(p,1)=\int_{\Bbb C}\frac 1{(1+\vert x_1\vert^2)^p}dm(x_1)
=\frac {\pi}{p-1}.\tag 3.11$$
Hence, following the notation in the introduction, $p_1=1$.  
The corresponding normalized measure, as in $(0.14)$, is 
$$\frac {l+1}{\pi}\frac 1{(1+\vert x_1\vert^2)^{1+1+l}}dm(x_1)\tag 3.12$$
The Fourier transform $u(y_1)$ is the multiple of the 
fundamental solution for 
$$(1+\Delta_{y_1})^{1+1+l}\tag 3.13$$
normalized so that $u(0)=1$, i.e.  
$$(1+\Delta )^{1+1+l}u(y_1)=\frac {l+1}{\pi}\delta_0(y_1).\tag 3.14$$

(2) For $N=2$ 
$$\int\frac 1{(1+\vert x_1\vert^2+2\vert x_2\vert^2+\vert x_2\vert^
4)^p}dm(x_1,x_2)=\tag 3.15$$
$$\int\frac 1{(\vert x_1\vert^2+(1+\vert x_2\vert^2)^2)^p}dm(x_1,
x_2)=\tag 3.16$$
$$\int\frac {dm(x_1')}{(1+\vert x'_1\vert^2)^p}\cdot\int\frac {dm
(x_2)}{(1+\vert x_2\vert^2)^{2p-2}}\tag 3.17$$
where we factored out the $(1+\vert x_2\vert^2)$ term, and made a 
change of variable in $x_1$.  Thus 
$$\Cal Z(p,2)=\frac {\pi^2}{(p-1)(2p-3)}.\tag 3.18$$
and $p_2=3/2$.  Thus the normalized measure, as in $(0.14)$, 
is 
$$\frac {2(l+1)(\frac 32+l)}{\pi^2}\frac 1{(\vert x_1\vert^2+(1+\vert 
x_2\vert^2)^2)^{1+\frac 32+l}}dm(x_1,x_2)\tag 3.19$$
The Fourier transform $u(y_1,y_2)$ is the multiple of the 
fundamental solution for 
$$(\Delta_{y_1}+(1+\Delta_{y_2})^2)^{1+\frac 32+l}\tag 3.20$$
normalized so that $u(0)=1$, i.e.  
$$(\Delta_{y_1}+(1+\Delta_{y_2})^2)^{1+\frac 32+l}u(y_1,y_2)=\frac {
2(l+1)(l+\frac 32)}{\pi^2}\delta_0(y_1,y_2).\tag 3.21$$

It seems unlikely that there is a simple factorization 
which allows one to directly evaluate the integrals when 
$2<N$.  However a miracle occurs.  There is a very 
simple torus action 
$$\Bbb T^2\times \{x=\sum_1^Nx_jz^j\}\to \{\sum_1^Nx_jz^j\}:((\lambda 
,\mu ),x)\to\sum_1^N\mu x_j\lambda^jz^j.\tag 3.22$$
It turns out that this is a Hamiltonian action, relative 
to a nonobvious symplectic form, with origins in the 
theory of loop groups, and $logdet(1+BB^{*})$ is a component 
of the momentum map.  Consequently one can apply the 
Duistermaat-Heckman theorem to evaluate the integrals.  
We will simply state the result; the details will appear 
in [Pi3].  

\proclaim{(3.23) Proposition} In general 
$$\int\frac 1{det(1+B_NB_N^{*})^p}dm(x_1,..,x_N)=\frac {\pi^N}{N!
(p-1)(p-(2-\frac 12))..(p-(2-\frac 1N))}.$$
Hence, in the notation of $(0.14)$, $p_N=2-\frac 1N$.  
\endproclaim

\flushpar(3.24) Remarks.  (1) In addition to having an 
explicit formula, this says that we have an incredibly 
simple recursion relation 
$$\Cal Z(p,N+1)=\frac {\pi}{(N+1)(p-(2-\frac 1N))}\Cal Z(p,N).\tag 3.25$$
There must be some direct proof of this.  

(2) Note that ignoring the $\pi^N$ factor, $(3.23)$ equals 
$$\prod_{k=1}^N(kp-(2k-1))^{-1}=\prod_k(1+(p-2)k)^{-1}=\frac {z^N}{
N!}\prod_{k=1}^N\frac 1{1+\frac zk}$$
$$=\frac {z^N\Gamma (z+1)}{\Gamma (z+N+1)}\tag 3.26$$
where $z=\frac 1{p-2}$.  This is related to the $\Gamma$ function, 
$$\Gamma (z)=\frac {e^{-\gamma z}}z\prod_{k=1}^{\infty}(\frac 1{1
+\frac zk})e^{z/k},\tag 3.27$$

\proclaim{(3.28) Theorem} The sequence of probability 
measures 
$$\frac 1{\pi^N}\prod_{k=1}^N(1+(l+1)k)\frac 1{det(1+B_N(x)B_N(x)^{
*})^{1+p_N+l}}dm(x_1,..,x_N)$$
converge to an invariant measure on $H^1=H^0/\Bbb C$, for each 
$-1<l$.  
\endproclaim

\flushpar(3.29) Remarks.  (a) This would be a beautifully 
simple result if these measures indexed by $N$ were 
coherent.  But unfortunately this is not the case 
(although when one calculates (using Maple, say), one 
finds that they are close to being coherent, even for 
small $N$).  

(a) It is clearly important to compute the 
exact $P_N$ distributions, or their Fourier transforms.
Another possibility is
to consider the values of $\theta$ at points 
of a configuration in $\Delta$, and try to compute the 
distribution for $(\theta (z_1),..\theta (z_n))$.  

Corresponding to this configuration we consider 
$$\theta =\sum\frac {\theta (z_j)}{(1-\bar {z}_jz)^2}\tag 3.30$$
This is the function that represents evaluation at the 
points.  The corresponding 
$$x=\sum\frac {\theta (z_j)}{\bar {z}_j(1-\bar {z}_jz)}+constant\tag 3.31$$
$$=\sum\theta (z_j)\frac z{1-\bar {z}_jz}\tag 3.32$$
Unfortunately it is not clear whether this yields any 
advantage in trying to understand the corresponding 
integrals involving $det(1+BB^{*})$ for this form of $x$.  

\bigskip

\flushpar Generalized Hankel Operators.

\smallskip

Fix $m,n>0$.  Define 
$$B_N(m,n;F)=B(m,n;(\sum_0^{N-1}F_jz^j)(dz)^{m+n})\tag 3.33$$
We will abbreviate this to $B_N$.  

It is clear that for sufficiently large $p$, 
$$\frac 1{det(1+B_NB_N^{*})^p}\in L^1(dm(F_0,..,F_{N-1})).\tag 3.34$$
Let $p_N$ denote the infimum of these $p$, and 
$$\Cal Z(m,n,;p)=\int\frac 1{det(1+B_NB_N^{*})^p}dm(F_0,..,F_{N-1}
)<\infty ,\quad p_N<p.\tag 3.35$$

With these additional parameters floating around, it 
seems less likely that a localization argument will apply 
to calculate $\Cal Z$ exactly.  However there is some remote 
possibility that there is some relation to Selberg 
integrals.  

\proclaim{(3.23) Question} Is there some obvious way to 
generalize our formula for $\Cal Z$ so that when $m=n=\frac 12$ 
$$\Cal Z(m,n,N;p)=\frac {\pi^N}{N!(p-1)(p-(2-\frac 12))..(p-(2-\frac 
1N))}.$$
\endproclaim

Regardless of whether there is an exact formula in this 
generality or not, it should always be the case that 
$$p_c=p_c(m,n)=\sup_Np_N<\infty .$$

\proclaim{(3.28) Conjecture} For $p_c<p$, the sequence of 
probability measures 
$$\frac 1{\Cal Z(m,n,N;p)}\frac 1{det(1+B_N(m,n;F)B_N(m,n;F)^{*})^
p}dm(P_NF)$$
converge to an invariant measure on $H^{m+n}$.  
\endproclaim

It is straightforward to generalize these conjectures in 
various ways.  For example using the decomposition 
$$H^{m_1}\cap L^2\otimes ..\otimes H^{m_r}\cap L^2\otimes H^{n_1}
\cap L^2\otimes ..\otimes H^{n_s}\cap L^2=H^{\sum (m_j+n_j)}\cap 
L^2\oplus ..,\tag 3.29$$
generalizing $(3.9)$, there is a Hilbert-Schmidt operator 
$$B(\vec {m},\vec {n};F)\in \Cal L_2((H^{m_1}\cap L^2\otimes ..\otimes 
H^{m_r}\cap L^2)^{*},H^{n_1}\cap L^2\otimes ..\otimes H^{n_s}\cap 
L^2)\tag 3.30$$
for $F\in H^{\sum (m_j+n_j)}\cap L^2$.  All of these generalizations 
involve equivariant operator-valued maps which are 
linear.

\bigskip

\flushpar Other Generalizations.

\smallskip

Recall from Lemma $(2.9)$ that
$$H^{1/2}\cap L^2\otimes H^{1/2}\cap L^2=H^1\cap L^2\oplus H^2\cap 
L^2\oplus ...$$
From this perspective the map
$$H^m\cap L^2\to H^{1/2}\cap L^2\otimes H^{1/2}\cap L^2$$
is a natural generalization of the classical 
Hankel map $x\to B(x)$, which corresponds to $m=1$.
To appreciate the special nature of these 
maps, we go back to the beginning.

The infinitesimal action of $sl(2,\Bbb C)$ on $H^m$ is given by 
$$A=\left(\matrix 0&0\\
1&0\endmatrix \right)\to -\frac {\partial}{\partial z}$$
$$A^{*}=\left(\matrix 0&1\\
0&0\endmatrix \right)\to z^2\frac {\partial}{\partial z}+2mz\tag A.1$$
$$H=\left(\matrix 1&0\\
0&-1\endmatrix \right)\to 2z\frac {\partial}{\partial z}+2m$$
The spectrum of the Hamiltonian is 
$$2m,2m+2,2m+4,..\tag A.2$$
From this one can read off that the action is essentially 
irreducible, except when $m$ is a nonpositive half-integer.
For $s=0,\frac 12,1,..$, and $m=-s$, there is a short exact 
sequence 
$$0\to \Bbb C^{2s+1}=span\{z^0(dz)^m,..,z^{2s}(dz)^m\}\to H^m\to 
H^{1+s}\to 0.\tag A.7$$

We are interested in the `bosonic cases' $s=0,1,..$ 
The case $s=0$ is classical Hankel  
$$H^0@>{\partial}>>H^1:x\to\theta\tag A.3$$
As we have emphasized, the Hankel operator is thought 
of most naturally in terms of $x$ (this is also the 
linearization of the mapping we will consider in \S 4).

Suppose $s=1$. In this case $()$ is realized as the map
$$0\to sl(2,\Bbb C)\to H^{-1}@>{dS}>>H^2\to 0\tag A.5$$
where
$$dS:v=(v_2z^3+..)\frac {\partial}{\partial z}\to Q=v^{\prime\prime
\prime}(dz)^2=(Q_2+Q_3z+..)(dz)^2,$$
$$Q_2=3\cdot 2v_2,Q_3=4\cdot 3v_3,..$$
(This is the linearization of the Schwarzian 
$$SL(2,\Bbb C)\to exp(H^{-1})@>{S}>>H^2\tag A.6$$
where we think of $exp(H^{-1})$ as formal holomorphic 
automorphisms of the disk, and which we pursue in \S 5).

In this case, instead of considering the multiplicative action of $
x\in H^0$ 
on $\Omega^{1/2}$, we consider the natural 
infinitesimal action of the vector field $v$. Because
$$\Cal L_v(z^n(dz)^{1/2})=(v(z)nz^{n-1}+\frac 12v'(z)z^n)(dz)^{1/
2}$$
$$=(z^n\sum_2^{\infty}\frac 12(2n+j+1)v_jz^j)(dz)^{1/2}$$
This leads to
$$B(v)=\frac 12\left(\matrix 4v_5&&&&\\
3v_4&2v_5&..&&\\
2v_3&v_4&0&.&\\
v_2&0&-v_4&-2v_5&.\\
0&-v_2&-2v_3&-3v_4&-4v_5&..\endmatrix \right)$$

At this point I do not have any feeling for the behavior 
of the critical exponents.

In the next case we should seek a map making the 
following diagram commute:
$$\matrix &&H^{-2}\\
&\swarrow&&\searrow\\
H^3&&\to&&H^{1/2}\otimes H^{1/2}\endmatrix $$
and so on.

\bigskip

\centerline{\S 4. Loop Groups and Conformal Invariance.}

\bigskip

As discussed in the introduction, the theory of loop 
groups yields further examples of conformally invariant 
measures in $H^1$.  In place of $(0.11)$ and $(0.12)$, we will 
consider a more complicated equivariant operator-valued 
function 
$$\theta_{+}\to W(x_{+}),$$
and the associated invariant function 
$$det(1+W(x_{+})W(x_{+})^{*})\tag 4.1$$
where $\theta_{+}\in H^1(\Delta ,\frak k^{\Bbb C})$ and $x_{+}=\int
\theta_{+}$.  This leads to a 
family of 
invariant measures on $H^1(\Delta ,\frak k^{\Bbb C})$ of the form 
$$\frac 1{\Cal Z}\frac 1{det(1+WW^{*})^{\frac 1m(2\dot {g}+l)}}dm
(\theta_{+}),\tag 4.2$$
where $l\ge 0$, and the other symbols will be explained 
below.  

Let $K$ denote a connected compact Lie group, and let $K^{\Bbb C}$ 
denote the complexification.  We are mainly interested in 
the case when the Lie algebra $\frak k$ is simple, e.g.  
$K=SU(n)$.  But we will want to compare with the 
abelian case, $K=U(1)$, so we will not assume simplicity 
at the outset.  

Fix a representation of $K$ in $U(N)$.  As in \S 3, there is a 
Hilbert space isomorphism 
$$L^2(S^1,\Bbb C^N)\to\Omega^0(\kappa^{1/2}\otimes \Bbb C^N):f(z)
\to f(z)(dz)^{1/2},\tag 4.3$$
and the Hardy polarization of $L^2(S^1,\Bbb C^N)$ is identified 
with the $SU(1,1)$-equivariant polarization 
$$\Omega^{1/2}\cap L^2\otimes \Bbb C^N=(H^{1/2}\cap L^2\otimes \Bbb C^
N)\oplus (H^{1/2}(\Delta^{*})\cap L^2\otimes \Bbb C^{^N}).\tag 4.4$$

As in chapter 6 of [PS], given $g\in LK^{\Bbb C}$, we obtain a 
multiplication 
operator $M_g$ on $\Omega^{1/2}\cap L^2\otimes \Bbb C^N$.
Relative to the Hardy polarization $(4.4)$,  
$$M_g=\left(\matrix A&B\\
C&D\endmatrix \right),\tag 4.5$$
where $A$ is referred to as the Toeplitz operator 
associated to $g$, and $B$ is the Hankel operator. If we 
introduce the standard 
ordered orthonormal basis 
$$..,z^2(dz)^{1/2},z^1(dz)^{1/2},z^0(dz)^{1/2};z^{-1}(dz)^{1/2},z^{
-2}(dz)^{1/2},..$$
then $M_g$, and its Toeplitz and Hankel operators, are 
represented by infinite block matrices,
$$M_g=\left(\matrix .&.\\
g_0&g_1&g_2&\\
.&g_0&g_1&\vert&g_2&g_3\\
.&g_{-1}&g_0&\vert&g_1&g_2&g_3&.\\
&-&-&-&-&-\\
.&g_{-2}&g_{-1}&\vert&g_0&g_1&g_2&.\\
.&g_{-3}&g_{-2}&\vert&g_{-1}&g_0&g_1&.\\
.&.&.&\vert&.&.&.\endmatrix \right)\tag $4.5'$ $$
which are constant along diagonals, 
where $g=\sum g_nz^n$ is the Fourier series of $g$ as a 
$\Cal L(\Bbb C^N)$-valued function. We also 
introduce the graph operators $Z=CA^{-1}$ and $W=A^{-1}B$; 
these arise when we consider the matrix LDU 
factorization 
$$\left(\matrix A&B\\
C&D\endmatrix \right)=\left(\matrix 1&0\\
Z&1\endmatrix \right)\left(\matrix A&0\\
0&D-ZB\endmatrix \right)\left(\matrix 1&W\\
0&1\endmatrix \right).\tag 4.6$$

Because the Hardy polarization identifies with a 
$SU(1,1)$-invariant polarization, the operators $A$, $B$, $C$, $D$, $
Z$ 
and $W$, are all equivariant 
operator-valued functions, in an appropriate sense; 
for example, the map to the Toeplitz operator
$$LK^{\Bbb C}\to \Cal L(H^{1/2}\cap L^2\otimes \Bbb C^N):g\to A(g
)\tag 4.7$$
intertwines the natural automorphic action of $PSU(1,1)$ 
and the action of $SU(1,1)$ by conjugation. We will mainly 
focus on $W$ below. 

It is a basic fact that $A=A(g)$ is invertible if and only 
if $g\in LK^{\Bbb C}$ has a Riemann-Hilbert factorization 
$$g=g_{-}\cdot g_0\cdot g_{+},\tag 4.8$$
where $g_{-}\in H^0(D^{*},\infty ;K^{\Bbb C},1),$ $g_0\in K^{\Bbb C}$, and $
g_{+}\in H^0(D,0;K^{\Bbb C},1)$, 
and in this case 
$$A((g_0g_{+})^{-1})=1_{N\times N}\quad or\quad g_0g_{+}=(A^{-1}(
\epsilon_1),..,A^{-1}(\epsilon_N))^{-1}\tag 4.9$$
where $1_{N\times N}$ denotes the identity $N\times N$ matrix, we are 
applying the operator $A$ column by column in the first 
equation, and $\epsilon_i$ denotes the standard basis for $\Bbb C^
N$ in 
the second equation.  Consequently there is a 
$PSU(1,1)\times K^{\Bbb C}$-equivariant map 
$$\{g\in LK^{\Bbb C}:det\vert A(g)\vert^2\ne 0\}\to H^1(D,\frak k^{
\Bbb C}):g\to\theta_{+}=g^{-1}_{+}\partial g_{+}.\tag 4.10$$

In the following Lemma, which summarizes well-known 
facts, $c(\cdot ,\cdot )$ denotes the Kac-Moody cocycle defined in 
Proposition $(6.6.4)$, page $88$ of [PS].  

\proclaim{(4.11) Lemma}For $g\in LK$, with Riemann-Hilbert 
factorization as in $(4.8)$, 
$$det(A(g)^{*}A(g))=det(1+WW^{*})^{-1},\quad W=A^{-1}B$$
$$=det(A(g)A(g)^{*})=det(1+Z^{*}Z)^{-1},\quad Z=CA^{-1}$$
$$=det\{A(g_{+})A(g_{+}^{*})A(g_{+}g_{+}^{*})^{-1}\}=c(g_{+},g_{+}^{
*}).$$
and 
$$Z(g)=Z(g_{-}),\quad W(g)=W(g_{+}).$$
These determinants are finite precisely when $g$ belongs 
to the Sobolev space $W^{1/2}$.
\endproclaim

\demo{Proof of (4.11)} Note that 
$$M_g=\left(\matrix A(g)&B(g)\\
C(g)&D(g)\endmatrix \right)=\left(\matrix A(g_{-})&0\\
C(g_{-})&D(g_{-})\endmatrix \right)\left(\matrix A(g_0g_{+})&B(g_
0g_{+})\\
0&D(g_0g_{+})\endmatrix \right)\tag 4.12$$
This implies that 
$$Z(g)=C(g)A(g)^{-1}=C(g_{-})A(g_0g_{+})\{A(g_{-})A(g_0g_{+})\}^{
-1}=Z(g_{-})\tag 4.13$$
and similarly $W(g)=W(g_{+})$.  Therefore 
$$1+Z^{*}Z=1+(C(g_{-})A(g_{-})^{-1})^{*}C(g_{-})A(g_{-})^{-1})$$
$$=(A(g_{-})^{*}A(g_{-}))^{-1}A(g_{-}^{*}g_{-})=A(g_{-})^{-1}A(g_{
-}^{*})^{-1}A(g_{-}^{*}g_{-}).\tag 4.14$$
with a similar expression involving $g_{+}$. 

By inspecting the matrix expression $(4.5')$, it is apparent 
that
$B$ and $C$ are Hilbert-Schmidt if and 
only if $g$ belongs to $W^{1/2}$. Since $g$ is automatically 
bounded (because it is $K$-valued) , this leads to the Sobolev 
condition.
\qed
\enddemo

Before discussing determinants and measures, 
we will first consider the 
operator
$$W(g_{+})=A(g_{+})^{-1}B(g_{+})=A(g_{+}^{-1})B(g_{+}).\tag 4.15$$
Ultimately we are interested in understanding this as
an equivariant holomorphic operator-valued function 
of $\theta_{+}$ in $(4.10)$, but we will first compute it in terms of 
$g_{+}$. We can picture $W$ as a block matrix, $W=(W_{i,-j})$, 
where each block $W_{i,-j}\in \Cal L(\Bbb C^N)$, 
and the index $i\ge 0$ denotes the row 
(starting from the bottom) and $j\ge 1$ denotes the column 
(starting from the left). If 
$$g_{+}=1+g_1z+g_2z^2+..\tag 4.16$$
where $g_j\in \Cal L(\Bbb C^N)$, then
$$g_{+}^{-1}=1+(-g_1)z+(-g_2+g_1^2)z^2+(-g_3+g_1g_2+g_2g_1-g_1^3)
z^3+..$$
$$\tag 4.17$$
and in general
$$(g_{+}^{-1})_n=\sum (-1)^lg_{i_1}..g_{i_l}\tag 4.18$$
where the sum is over all $positive$ multi-indices 
$I=(i_1,..,i_l)$ of order $n$, i.e. $i_m>0$ and 
$i_1+..+i_l=n$. By $(4.15)$,
$$W_{i,-j}=1\cdot g_{i+j}+(g_{+}^{-1})_1g_{i+j-1}+..+(g_{+}^{-1})_
ig_j$$
$$=\sum_{n=j}^{i+j}\sum (-1)^lg_{i_1}..g_{i_l}g_n\tag 4.19$$
where given $n$, the sum is over all positive multi-indices 
of order $i+j-n$. This can also be written as
$$W_{i,-j}=\sum (-1)^{l+1}g_{i_1}..g_{i_l}\tag 4.20$$
where the sum is now over all positive multi-indices of 
order $i+j$ satisfying $i_l\ge j$. Thus, in terms of the 
representation dependent expression $(4.16)$ for $g_{+}$, 
$W$ has the form
$$\left(\matrix ...&...\\
g_3-g_1g_2-g_2g_1+g_1^3&g_4-g_1g_3-g_2g_2+g_1^2g_2&...&\\
g_2-g_1g_1&g_3-g_1g_2&g_4-g_1g_3\\
g_1&g_2&g_3&..\endmatrix \right)\tag 4.21$$

Now write
$$\theta_{+}=(\theta_1+\theta_2z+..)dz\in H^1(D,\frak g),\tag 4.22$$
where $\theta_i\in \frak g$. Since $g_{+}$ is the solution of the integral 
equation
$$g_{+}(z)=1+\int_0^zg_{+}(w)\theta (w),\quad g_{+}(0)=1,\tag 4.23$$
it can be expressed in terms of iterated integrals:
$$g_{+}(z)=1+\int\theta +\int \{(\int\theta )\theta \}+\int \{(\int
\theta )\int$$
$$=1+g^{(1)}(\theta )+g^{(2)}(\theta )+..\tag 4.24$$
where 
$$g^{(n)}(\theta )=\int g^{(n-1)}(\theta )\theta =\sum\frac 1{i_1}
..\frac 1{i_1+..i_n}\theta_{i_1}..\theta_{i_n}z^{\vert I\vert}\tag 4.25$$
and the sum is over all positive multiindices 
$I=(i_1,..,i_n)$. 

Given a positive multi-index $I=(i_1,i_2,..,i_l)$, define 
$$c(I)=\frac 1{i_1}\frac 1{i_1+i_2}..\frac 1{i_1+..+i_l}.\tag 4.26$$
Observe that there is a bijective correspondence between 
positive multi-indices $I$ of order $n$ and subsets of 
$S\subset \{1,..,n-1\}$:  A multi-index $I$ induces a strictly 
increasing sequence 
$$\lambda_1=i_1<\lambda_2=i_1+i_2<..<\lambda_l=i_1+..+i_l=n\tag 4.27$$
which is uniquely determined by the complement 
$$S=\{1,..,n\}\setminus \{\lambda_1,..,\lambda_l\}\subset \{1,..,
n-1\}.\tag 4.28$$
In terms of $S$, the $i_j$ are of the form 
$$i_j=1+\vert S_j\vert ,\tag 4.29$$
where $S_j$ is the $j$th connected component of $S$, and two 
integers are connected if they are adjacent.  We can 
then write 
$$c(I)=\frac {\prod_S\lambda}{n!}\tag 4.30$$
We have 
$$g_n=\sum c(I)\theta_I=\sum c(I)\theta_{i_1}..\theta_{i_l}\tag 4.31$$
where the sum is over all positive multi-indices of 
order $n$. Plugging this into $(4.20)$ implies the following 
formula.

\proclaim{ (4.32) Proposition}As a function of $\theta_{+}\in H^1
(D,\frak g)$,
for $i\ge 0$ and $j\ge 1$,
$$W_{i,-j}=\sum C(I)\theta_{i_1}..\theta_{i_l},$$
where $I$ ranges over all positive multi-indices of order 
$i+j$, and
$$C(I)=\sum (-1)^{l+1}c(I_1)..c(I_l)$$
$$=\sum (-1)^{l+1}\frac {\prod_{S_1}\lambda_1}{n!}..\frac {\prod_{
S_l}\lambda_l}{n!}$$
where the sum is over all ways of representing $I$ as a 
tuple $(I_1,..,I_l)$ with $\vert I_l\vert\ge j$.  
\endproclaim

Thus $W(\theta_{+})$ has the form
$$\left(\matrix ..&..\\
\frac 1{3!}(2\theta_3-\theta_1\theta_2-2\theta_2\theta_1+\theta_1^
3)&..&..&\\
\frac 12(\theta_2-\theta_1\theta_1)&\frac 1{3!}(2\theta_3-[\theta_
1,\theta_2]-2\theta_1^3)&..\\
\theta_1&\frac 12(\theta_2+\theta_1\theta_1)&\frac 1{3!}(2\theta_
3+2\theta_1\theta_2+\theta_2\theta_1+\theta_1^3)&..\endmatrix \right
)$$

Now that we have a formula, the natural question is 
whether we can understand the terms of the Taylor 
series, in terms of the representation theory of 
$PSU(1,1)$. The Taylor series is
$$W=W^{(1)}+W^{(2)}+W^{(3)}+..$$
where $W^{(n)}$ is a homogeneous symmetric function of the 
$\theta_i$.

For example $W^{(1)}$ is the natural inclusion
$$H^1(D;\frak g)\to H^{1/2}\otimes H^{1/2}\otimes \Cal L(\Bbb C^N
)=(H^1+H^2+..)\otimes \Cal L(\Bbb C^N).$$
The second derivative $W^{(2)}$ can be viewed as a linear function
$$S^2(H^1(D;\frak g))\to (H^1+..)\otimes \Cal L(\Bbb C^N)$$
and this time the function is not determined by 
equivariance. The calculations indicate that there should 
something very interesting to be said about these 
functions. But what is it?

We now consider the determinants in $(4.9)$. These reduce 
to familiar invariant functions in special cases.  
In the abelian case, $K=U(1)$, 
$$A(g_{+})A(g_{+}^{*})A(g_{+}g_{+}^{*})^{-1}=$$
$$=A(g_{+})A(g_{+}^{*})A(g^{*}_{+}g_{+})^{-1}=$$
$$=A(g_{+})A(g_{+}^{*})A(g_{+})^{-1}A(g_{+}^{*})^{-1}=$$
$$=e^{A(x_{+})}e^{A(x_{+}^{*})}e^{-A(x_{+})}e^{-A(x_{+}^{*})}.\tag 4.33$$
This last expression is a commutator, and consequently 
the determinant can be calculated using the Helton-Howe 
formula ([HH]):  
$$c(g_{+},g_{+}^{*})=det([e^{A(x_{+})},e^{A(x_{+}^{*})}])=$$
$$e^{tr[A(x_{+}),A(x_{+}^{*})]}=e^{\sum_{n=1}^{\infty}n\vert x_n\vert^
2}=e^{\vert\theta_{+}\vert^2_{H^1(\Delta )}}.\tag 4.34$$

Now suppose that $K=SU(2)$.  If $g\in LK$ has a diagonal 
form, $g=\left(\matrix a&0\\
0&\bar {a}\endmatrix \right)$, then the abelian case implies 
$$det\vert A(g)\vert^2=exp(\vert\theta_{+}\vert^2_{H^1\cap L^2}),\tag 4.35$$
where $\theta_{+}$ is now a matrix of the special form 
$\theta_{+}=\left(\matrix \alpha&0\\
0&-\alpha\endmatrix \right)$, $\alpha\in H^1(\Delta )$.  

At the opposite extreme is the nilpotent case:  
$$\theta_{+}=\partial x_{+}=\left(\matrix 0&\theta\\
0&0\endmatrix \right),\quad\theta\in H^1(\Delta ),\tag 4.36$$
where 
$$g_{+}=\left(\matrix 1&\int^z_0\theta\\
0&1\endmatrix \right)=\left(\matrix 1&x\\
0&1\endmatrix \right).\tag 4.37$$

\proclaim{(4.38)Lemma} Relative to the defining 
representation, 
$$W(g_{+})=B(g_{+})=\left(\matrix .\\
0&x_n&.&\\
0&0&.\\
.&.\\
0&x_2\\
0&0&.&&&.&.\\
0&x_1&0&x_2&&.&&0&x_n&.\\
0&0&0&0&..&&&0&0&.\endmatrix \right).$$
$$\roman w\roman h\roman i\roman c\roman h,\ \roman m\roman o\roman d
\roman u\roman l\roman o\ \roman v\roman a\roman n\roman i\roman s
\roman h\roman i\roman n\roman g\ \roman c\roman o\roman l\roman u
\roman m\roman n\roman s\ \roman a\roman n\roman d\ \roman r\roman o
\roman w\roman s,\ \roman i\roman s\ \roman t\roman h\roman e\ \roman s
\roman a\roman m\roman e\ $$
as $B(x)$ in $(0.11)$.  
\endproclaim

\demo{Proof} One can either read this off from 
$(4.32)$, or argue directly as follows. Suppose that $\left(\matrix 
f_1\\
f_2\endmatrix \right)\in H^{-}$.  
$$A(g_{+})^{-1}B(g_{+})\left(\matrix f_1\\
f_2\endmatrix \right)=P_{+}g_{+}^{-1}P_{+}g_{+}\left(\matrix f_1\\
f_2\endmatrix \right)=$$
$$=P_{+}g_{+}^{-1}\left(\matrix f_1+P_{+}(xf_2)\\
f_2\endmatrix \right)=P_{+}\left(\matrix f_1+P_{+}(xf_2)-xf_2\\
f_2\endmatrix \right)=$$
$$\left(\matrix P_{+}(xf_2)\\
0\endmatrix \right)=B(g_{+})\left(\matrix f_1\\
f_2\endmatrix \right).\tag 4.39$$
The matrix for $B$ is 
linear as a function of $g_{+}$, and it vanishes on the identity. 
Since 
$g_{+}=1+xe_2^{*}\otimes e_1$, this implies that $B(g_{+})$ is $B
(x)$, modulo 
vanishing columns and rows.  \qed
\enddemo

These calculations show that in the abelian case the 
measure $(4.2)$ is a Gaussian (the number $\dot {g}=0$ in the 
abelian case, and we must require $l>0$), and if we 
consider nilpotent $\theta_{+}$, when $\frak g=sl_2$, 
then $(4.2)$ is to be interpreted as in \S 3.  

We now suppose that $\frak k$ is a simple Lie algebra.  Let $\dot {
g}$ 
denote the dual Coxeter number of $\frak k$; if $K=SU(n)$, then 
$\dot {g}=n$.  The number $m$ in $(4.2)$ depends upon $N$:  for 
$X,Y\in \frak k^{\Bbb C}$, 
$$tr_{\Bbb C^N}(XY)=\frac 1m\kappa (X,Y),\tag 4.40$$
where $\kappa$ denotes the Killing form for $\frak k^{\Bbb C}$ (if $
\Bbb C^N$ is the 
adjoint representation, then $m=\dot {g}$, and if $K=SU(n)$, and 
$N=n$, then $m=1$).  

\proclaim{ (4.41) Theorem} Fix $l\ge 0$.  (a) Fix $M$.  Given 
$\theta_{+}=P_M\theta_{+}=\sum^M\theta_iz^{-i}dz^{-1}$, let $g\in 
LK$ satisfy $g=g_{-}$ 
modulo $H^0(D,G)$.  The measure 
$$\frac 1{det(1+W(\theta_{+})W(\theta_{+})^{*})^{\frac 1m(2\dot {
g}+l)}}dm(P_M\theta_{+})\tag 4.42$$
is finite.  

(b) Let $\Cal Z_M$ denote the total mass.  The limit of the 
probability measures 
$$\frac 1{\Cal Z_M}\frac 1{det(1+WW^{*})^{\frac 1m(2\dot {g}+l)}}
dm(P_M\theta_{+})\tag 4.43$$
exists (in a weak sense), and defines a $PSU(1,1)$-invariant 
measure on $H^1(\Delta ,\frak k^{\Bbb C})$.  Consequently, for each matrix 
coefficient of $\theta_{+}$, we obtain an invariant measure on $H^
1$.  
\endproclaim

This will appear in [Pi4].  

Because these formulas involve limits, they seem to be 
of limited utility, beyond heuristically explaining why 
the measures are $PSU(1,1)$-invariant. For example
it seems plausible that the measure in $(0.14)$ is the 
invariant measure that corresponds to the highest root,
but this is uncertain.  These measures are so natural,
from a group point of view, that one has to believe 
there is some direct way of computing them.

\bigskip

\centerline{\S 5. Universal Teichmuller Space.}

\bigskip

In this section we will consider another invariant 
function, this time associated with diffeomorphisms of 
the circle.  The function will have the same form as in 
$(4.1)$, but it is far more complicated.  

The first object of this section is to recall the Bers 
embedding of universal Teichmuller space into $H^2(\Delta )$.  
This embedding is equivariant with respect to the 
natural actions of $PSU(1,1)$ on these spaces.  The 
hypothetical invariant measures which we will discuss 
should be supported on the image of this embedding.  

We will use the following notations.  The invariant norm 
$(0.24)$ induces a natural $PSU(1,1)$-invariant $L^{\infty}$ norm on 
quadratic differentials:  
$$\vert Q(z)(dz)^2\vert_{L^{\infty}}=\sup_{\Delta}\{(1-\vert z\vert^
2)^2\vert Q(z)\vert \}.\tag 5.1$$
To gain a feeling for this norm, note that a quadratic 
differential which extends meromorphically to a 
neighborhood of $\bar{\Delta}$, and which is bounded in this 
$PSU$-invariant sense, can have at most double poles 
around the boundary, $S^1$.  We will denote this Banach 
space by $H^2\cap L^{\infty}$.  

We will denote the space of normalized univalent 
holomorphic functions on $\Delta$ by $\Cal S$:  $u\in \Cal S$ means that $
u$ is 
a $1-1$ holomorphic function on $\Delta$ having a Taylor series 
expansion of the form 
$$u(z)=z(1+\sum_1^{\infty}u_nz^n).\tag 5.2$$
The space of orientation-preserving homeomorphisms of 
$S^1$ will be denoted by $\Cal D$, and we will use subscripts to 
denote further restrictions on the degree of smoothness.  
For our purposes the most important example is $\Cal D_{qs}$, the 
group of quasi-symmetric homeomorphisms.  The 
condition that $\sigma$ is quasi-symmetric can be expressed in 
several equivalent ways:  (1) $\sigma$ is the restriction to $S^1$ 
of a quasi-conformal automorphism of $\hat {\Bbb C}$ (or $\Delta$) mapping 
$S^1$ into itself; (2) $\sigma$ satisfies the Beurling-Ahlfors 
criterion (see chapter 16 of [GL], especially $(16.1)$); and (3) 
$\sigma$ stabilizes $W^{1/2}(S^1)$ (see [NS], especially Theorem 3.1).  

The Bers embedding depends upon the following theorem 
of Bers (see [L], page 100).  

\proclaim{(5.3)Theorem}If $\sigma\in \Cal D_{qs}$, then there is a unique 
factorization (where the multiplication is composition of 
maps) 
$$\sigma =l\circ diag\circ u,\tag 5.4$$
where $diag$ is multiplication by a constant $\lambda\in \Bbb C^{
*}_{\le 1}$, 
$$u=z(1+\sum_1^{\infty}u_nz^n),\quad\frac 1{l^{-1}(\frac 1w)}=w(1
+\sum_1^{\infty}l_nw^n),\tag 5.5$$
$u$ is univalent in $\Delta$ and admits a quasi-conformal 
extension to $\hat {\Bbb C}$, $l^{-1}$ is univalent in $\Delta^{*}$ and admits a 
quasi-conformal extension to $\hat {\Bbb C}$, and the compatibility 
condition 
$$\lambda u(S^1)=l^{-1}(S^1)\tag 5.6$$
is satisfied.  
\endproclaim

The complement of $l^{-1}(\Delta^{*})$, $\lambda u(\Delta )$, has unit transfinite 
diameter, $\rho (\lambda u(\Delta ))=1$ (see \S 16.2 of [Hi] for the original 
definition of $\rho$, Fekete's theorem 16.2.2, page 270 for the 
basic characterization of $\rho$, and see the second paragraph 
of page 347 for this specific fact), hence $\frac 1{\vert\lambda\vert}
=\rho (u(\Delta ))$.  
An immediate consequence of this is that $u$ (or $l$) 
determine $\sigma$ up to a phase factor.  In particular we can 
think of $u$ as a parameter for the homogeneous space 
$Rot\backslash \Cal D_{qs}$.  

This leads to the following diagram:  
$$\matrix Rot\backslash \Cal D_{qs}&\leftrightarrow&\Cal S_{qc}&\hookrightarrow&
\Cal S&&\hookrightarrow&&\Cal S_{formal}\\
\downarrow&&\downarrow&&\downarrow&&&&\downarrow S\\
PSU\backslash \Cal D_{qs}&\leftrightarrow&\Cal T&\hookrightarrow&
\Cal U&\hookrightarrow&H^2(\Delta )&\hookrightarrow&H^2(\Delta )_{
formal}\endmatrix \tag 5.7$$
Here $\Cal S$ denotes the space of all univalent holomorphic 
functions $u$ on $\Delta$ of the form $(1.11)$, $\Cal S_{qc}$ denotes the 
subspace of all such $u$ which have quasi-conformal 
extensions to $\hat {\Bbb C}$, $\Cal S_{formal}$ consists of formal power series 
$u$ as in $(5.5)$, $S$ is the Schwarzian derivative, and $\Cal T$, the 
universal Teichmuller space, and $\Cal U$ are the $S$-images of 
$\Cal S_{qc}$ and $\Cal S$, respectively.  

The action of $PSU\subset \Cal D_{qs}$ from the right on $\Cal S_{
qc}$ extends 
to an action on $\Cal S$ which is given as follows.  Suppose 
that $u\in \Cal S$ and $g\in PSU$.  There exists a unique 
$G\in PSL(2,\Bbb C)$ such that (i) $G(0)=u(g(0))$, (ii) 
$G'(0)=(u\circ g)'(0)$ and (iii) $G(\infty )=\infty$.  The right action by 
$g$ is given by 
$$g:u\to G^{-1}\circ u\circ g.\tag 5.8$$
The cocycle condition 
$$S(f\circ g)=S(f)\vert_{g(z)}g'(z)^2+S(g),\tag 5.9$$
and the vanishing of $S$ on linear fractional maps implies 
that $S(G^{-1}\circ u\circ g)=S(u\circ g)$.  Hence the down arrows 
above are $PSU$-equivariant (excluding the rightmost 
arrow, because the action of $PSU$ does not extend to the 
formal completions).  Thus both $\Cal T$ and $\Cal U$ are invariant 
for the natural (right) pullback action of $PSU$ on $H^2(\Delta )$.  

We have implicitly observed here that $S$ extends to an 
algebraic map of the formal completions (see below for 
some speculation about the algebraic properties of this 
mapping).  

It is known that 
$$\{q\in H^2(\Delta )\cap L^{\infty}:\vert q\vert_{L^{\infty}}<2\}
\subset \Cal T\subset \Cal U\subset \{q\in H^2(\Delta )\cap L^{\infty}
:\vert q\vert_{L^{\infty}}<6\}\tag 5.10$$
Also $\Cal T$ is a bounded open subset of the Banach space 
$H^2\cap L^{\infty}$, so in this sense it is a homogeneous bounded 
domain.  A lot of work has gone into understanding the 
nature of the boundary of $\Cal T$, and its subspaces $\Cal T_{\Sigma}$, as $
\Sigma$ 
varies over all hyperbolic type Riemann surfaces.  The 
state of the art in understanding these spaces is 
described in [GL].  

To understand the Bers embedding, or the Schwarzian 
derivative, algebraically, we introduce the following 
mappings:  
$$\Cal S@>{c}>>H^0(\Delta ):u\to log(u'),\tag 5.11$$
and 
$$\Cal S@>{N}>>H^1(\Delta ):u\to log(u')'.\tag 5.12$$
These mappings satisfy the cocycle conditions 
$$c(f\circ g)\vert_z=c(f)\vert_{g(z)}+c(g)\vert_z\tag 5.13$$
$$N(f\circ g)=N(f)\vert_{g(z)}g'(z)+N(g)\vert_z\tag 5.14$$
and they have the property that they are equivariant 
with respect to the actions of $PSU(1,1)$ on $\Cal S$ and the 
natural actions on the $H^{*}(\Delta )$ spaces.  

We have 
$$u=z(1+\sum_1^{\infty}u_nz^n)\quad@>{S}>>\quad Q=\sum_0^{\infty}
Q_{n+2}z^n(dz)^2\tag 5.15$$
$$Q=ln(u')^{\prime\prime}-\frac 12ln(u')^{\prime 2}=N(u)'-\frac 1
2N(u)^2.\tag 5.16$$
Here are some sample calculations:  
$$c(u)=ln(u')=p_1(u)z+p_2(u)z^2+..,\quad$$
$$p_1=2u_1,p_2=3u_2-2u_1^2,p_3=4u_3-6u_1u_2+\frac 83u_1^3$$
$$p_4=5u_4-\frac 12(2\cdot 2u_14u_3+3u_23u_2)+\frac 13(3\cdot 2u_
12u_13u_2)-\frac 14(16u_1^4),..$$
$$N(u)=(lnu')'=p_1+2p_2z+3p_3z^2+4p_4z^3+..$$
$$Q=S(u)=N(u)'-\frac 12N(u)^2$$
$$=(2p_2+3\cdot 2p_3z+4\cdot 3p_4z^2+..)$$
$$-\frac 12(p_1^2+2p_12p_2z+[2p_13p_3++2p_22p_2]z^2+$$
$$[p_14p_4+2p_23p_3+3p_32p_2+4p_4p_1]z^3+..)$$
$$Q_2=2p_2-\frac 12p_1^2=6(u_2-u_1^2),$$
$$Q_3=6p_3-2p_1p_2=6(4u_3-6u_1u_2+\frac 83u_1^3)-2(2u_1)(2u_2-2u_
1^2)$$
$$=24u_3-14u_1u_2+20u_1^3.$$
$$Q_4=12p_4-3p_1p_3-2p_2^2$$
$$=12(5u_4-8u_1u_3+\frac 92u_2^2+12u_1^2u_2-4u_1^4)-$$
$$3(2u_1)(4u_3-6u_1u_2+\frac 83u_1^3)-2(3u_2-2u_1^2)^2$$

If we set $u_1=0$, then it seems fairly certain that there 
should be some tractable combinatorial description of the 
coefficients of $u$ in terms of $Q$:
$$u_2=\frac 16Q_2,\quad u_3=\frac 1{24}Q_3$$
$$u_4=\frac 1{5!}(2Q_4-2^23Q_2^2),..$$

Our goal now is to develop this section in analogy with 
the loop group case.

As we noted in \S 1B, the representation 
$$SU(1,1)\times H^{1/2}(\Delta )\to H^{1/2}(\Delta )\tag 5.18$$
is special.  If we consider the odd spin structure on $S^1$, 
then we obtain an action 
$$\Cal D^{(2)}\times\Omega^{1/2}_{odd}\to\Omega^{1/2}_{odd}.\tag 5.19$$
In concrete terms, if we introduce the double cover 
$$p:(S^1)^{(2)}\to S^1:\zeta\to z=\zeta^2,\tag 5.20$$
then 
$$\Cal D^{(2)}=\{\sigma\in Diff((S^1)^{(2)}):\sigma (-\zeta )=\sigma 
(\zeta )\}\tag 5.21$$
$$\Omega^{1/2}_{odd}=\{f(\zeta )d\zeta :f(-\zeta )=-f(\zeta )\}\tag 5.22$$
and $\Cal D^{(2)}$ acts in the natural way.  

There is a $SU(1,1)$-invariant polarization 
$$\Omega^{1/2}_{odd}=H^{1/2}(\Delta^{*})\oplus H^{1/2}(\Delta )\tag 5.23$$
of Hilbert spaces.  Given $\sigma\in \Cal D^{(2)}$, we write 
$$\sigma =\left(\matrix A&B\\
C&D\endmatrix \right)\tag 5.24$$
in analogy with the loop group case (hence we refer to 
$A$ as the Toeplitz operator corresponding to $\sigma$, etc.). 

The analogue of $(4.11)$ holds, where the Riemann-Hilbert 
factorization is replaced by the Bers factorization. 
Whereas it seems not so enlightening to examine a 
matrix representation for $(5.24)$, it is quite interesting 
to stare that the matrix representation for $u^{-1}$, relative 
to the standard basis 

$$..,z^2(dz)^{1/2},..,z^{-2}(dz)^{1/2},..\tag 5.25$$
It is a unipotent upper triangular matrix which the 
$(i,j)$-entry is homogeneous of degree $j-i$, where $u_l$ has 
degree $l$. It is very easy to generate this matrix, but 
not so easy to display it:
$$\left(\matrix .&3u_3+u_1u_2&\frac 52u_4-2u_1u_3-\frac 98u_2^2+\frac 
94u_1u_2-\frac 58u_1^4&&\\
3u_1&\frac 52u_2+\frac 12u_1^2&2u_3-\frac 32u_1u_2+\frac 12u_1^3&
\vert&\frac 32u_4-3u_1u_3-\frac {13}8u_1^4+\frac {19}4u_1^2u_2-\frac {
13}8u_2^2&\\
1&2u_1&\frac 32u_2-\frac 12u_1^2&\vert&u_3-2u_1u_2+u_1^3&\\
&1&(\cdot )_{(1,0)}=u_1&\vert&\frac 12u_2-\frac 12u_1^2&\\
&&(\cdot )_{(0,0)}=1&\vert&(\cdot )_{(0,-1)}=0&-\frac 12u_2+\frac 
12u_1^2&\\
-&-&-&-&-&-\\
&&&\vert&1&-u_1&..\\
&&&&0&1&-2u_1\endmatrix \right)$$

The key point is that there is a precise analogue of the 
linear coordinate for $g_{+}$, $\theta_{+}=g_{+}^{-1}\partial g_{
+}$ (see $(4.10)$).  In 
terms of the Bers embedding theorem, the analogue is 
the map 
$$\sigma\to Q=S(u)(dz)^2\in \Cal T\subset H^2\cap L^{\infty}\tag 5.26$$
Just as we studied the equivariant mapping
$$\theta_{+}\to W(g_{+})=A(g_{+}^{-1})B(g_{+})$$
we would like to be able to say something about the 
combinatorics of the mapping 
$$Q\to W(u)=A(u^{-1})B(u)$$
Here we will only say that it has the homogeneous form
$$\left(\matrix .&&\\
\frac 1{12}Q_4+cQ_2^2&..&&\\
-\frac 1{12}Q_3&-\frac 1{36}Q_4+c'Q^2_2&0&..&&&\\
\frac 13Q_2&0&\frac 1{36}Q_4+c^{\prime\prime}Q_2^2&&\\
0&-\frac 13Q_2&\frac 1{12}Q_3&-\frac 1{12}Q_4+dQ_2^2\endmatrix \right
)$$

The major new complications are that (1) whereas we have at
least a crude grasp of the correspondence between $\theta_{+}$ and
the analytic properties of $W(g_{+})$, this is quite 
mysterious in the case of
$Q$ and $W(u)$, and (2) whereas it was reasonable to use 
`infinite dimensional Lebesgue measure' as a heuristic 
background measure for $\theta_{+}$, we now need a background 
measure for $Q$ which is presumably supported on $\Cal T$.

We want to consider the critical exponents for the 
expression 
$$\frac 1{det(1+WW^{*})^p}dV(Q)\tag 5.28$$
over a sequence of finite dimensional space $\Cal T_N$  which 
tend to $\Cal T$ in the limit. One suspects that there must be 
many different interesting possibilities; the hard question is 
whether any of them lead to tractable calculations. At 
this point I will just mention one obvious possibility.

By $(5.7)$ we can identify $\Cal T$ with 
domains in $\Bbb C$ with Jordan curve boundaries, where two 
regions are identified if one can be obtained from 
another by a linear fractional transformation of $\hat {\Bbb C}$.
We let $\Cal T_N$ denote the subset corresponding to $N$-gons.
 
Given an $N$-gon, let $z_1$$,z_2,$$..,z_N$ denote the vertices 
listed in the counterclockwise direction. We can write
the interior angle of the polygon at $z_j$ in the form $\alpha_j\pi$, 
$0<\alpha_j<2$. We also let $\beta_j=1-\alpha_j$.  

In this case we can say exactly what the $u$ and $l$ are 
(we do not quite know $\lambda$), using the Schwarz-Chrisstoffel 
transform theory, namely we have 
$$du=\prod_1^n(1-\frac z{z_j})^{-\beta_j}=1+\sum\frac {\beta_j}{z_
j}z+..\quad\implies\quad u=z(1+u_1z+..)\tag 5.29$$
$$dl^{-1}=\prod_1^n(1-\frac {z'_j}z)^{\beta_j}=1+\frac {\sum\beta_
jz'_j}z+..\quad\implies\quad l^{-1}=z+b_0+b_1z^{-1}+..\tag 5.30$$
This enables one to do calculations. But the main thing 
that is lacking is a feeling for the background geometry,
which should approximate the invariant (Weil-Petersson)
geometry on universal Teichmuller space.

\bigskip

\centerline{\S 6. Comments and Questions}

\bigskip

In \S 1 we established the existence of transitive 
measures, hence of measures having finite dimensional 
support, for $H^m$, when $m>0$ and rational.  

$(6.1)$.  Do there exist measures having finite dimensional 
support when $m$ is irrational?  (We know there do not 
exist transitive measures in this case, so one suspects 
the answer is `no').  

$(6.2)$.  Suppose that $m>0$ is rational.  Is the convolution 
algebra generated by transitive measures dense in the 
space of all invariant (not necessarily positive) measures 
on $H^m(\Delta )$?  If yes, we could further ask if there is some 
canonical way to decompose an invariant measure of 
finite dimensional support.  If no, we could allow 
multiplication of such measures and ask the same 
questions.  

The results in [K], summarized in $(1.15)$, go a long 
way toward describing the support of the $G$-ergodic 
Gaussian measure corresponding to the Hilbert space 
$H^1\cap L^2$.  

$(6.3a)$.  Can this be refined to say in an explicit 
way how the dynamics of $G$ is chaotic?  

$(6.3b)$. Is there something interesting to be said about
$\{\Bbb Re(x)=0\}$. This is a curve that starts at zero and goes 
out to $\infty$ in a random. This is the opposite of chordal 
SLE, where one starts at infinity and goes towards the 
origin. We are also looking at the whole curve at once; 
there is not a time parameter. We could ask this 
question for any of the invariant distributions on 
functions.

$(6.)$. Are there examples of invariant distributions for which 
one can compute the distribution of zeroes? 
Does the distribution of zeroes determine 
the distribution on $H^m$, or 
at the other extreme, possibly there are conditions 
under which it is independent of the $H^m$-distribution.

The following table is a rather feeble attempt to 
indicate some of what we do and do not know in an 
explicit way, regarding invariant measures for $H^1=H^0/\Bbb C$ 
($y$ is a variable dual to $x\in H^0/\Bbb C$):  
 
$$\matrix \text{\rm $\underline {\text{\rm Measure}}$}&&\underline {\text{\rm Fourier Transform}}\\
\text{\rm a transitive measure}&&?\\
&&\\
\nu^{(1)}_t=\frac 1{\Cal Z}e^{-\frac 1{2t}\vert x\vert^2}d\lambda 
(x)&&e^{-\frac t2\vert y\vert^2}\\
(\nu^{(\lambda )}_s\otimes\nu_t^{(1-\lambda )})^{*N}&&\frac 1{det
(1+stB(\lambda ,1-\lambda ;\cdot )B(\lambda ,1-\lambda ;\cdot )^{
*})^N}\\
\nu^{(\lambda_1)}\otimes\nu^{(\lambda_2)}\otimes\nu^{(\lambda_3)}
\otimes ..&&?\\
\frac 1{\Cal Z}\frac 1{det(1+B(x)B(x)^{*})^p}d\lambda (x)&&u(y),\quad 
det(1+B(\frac {\partial}{\partial y})B(\frac {\partial}{\partial 
y})^{*})^pu=c\delta_0\\
\text{\rm matrix coeffs of }\frac 1{\Cal Z}\frac {d\lambda (x)}{d
et(1+WW^{*})^{2\dot {g}+l}}&&?\endmatrix $$

$(6.4)$.  In the second to last line, we have indicated that 
the Fourier transform is heuristically a `fundamental 
solution' for a constant coefficient operator (involving 
infinitely many variables).  Can one make sense of this?  

$(6.5)$.  If we could fill in the question mark on the 
fourth line, the table should continue in some interesting 
way.  

The Lebesgue type measures which we have come across 
in this paper have the property that their Fourier 
transforms are continuous on $(H^m\cap L^2)^{*}$.  It appears that 
the Fourier transform for an invariant measure of finite 
dimensional support is not continuous on $H^m\cap L^2$.  

$(6.6)$ Does this continuity property capture what should 
be meant by having `Lebesgue type support'?  

In this paper we have encountered a lot of invariant 
measures on $H^1$ which have $x_1$-density of the form 
$\Cal Z^{-1}(1+\vert x_1\vert^2)^{-p}$.  
The measures in \S 4, involving loops into $SU(2)$, or more 
generally a compact Lie group, have matrix coefficients 
which have this property.  

$(6.7)$ If one considers loops into $S^2$, one expects to 
encounter a measure having a matrix coefficient with 
$x_1$-density 
$$\frac 1{\Cal Z}\frac 1{(1+\vert x_1\vert^2)^{3/2}}F(\frac 1{1+\vert 
x_1\vert^2})$$
where $F(\rho )$ is the function 
$$F(\rho )=\int_0^{\infty}\frac {\rho}{(\rho +\frac {(x-1)^2}x)^{
3/2}}dx$$
How does this conjectural measure fit into our scheme?  

$(6.8)$ Referring to \S 4, is $\frac 1{det(1+WW^{*})}$ a positive definite 
function on $H^1\cap L^2\otimes \frak k^{\Bbb C}$?  If so, what is the 
corresponding measure?  

$(6.9)$ Consider the finite Hankel matrix $B_n(x)$.  Can one 
find a reasonable formula for the distribution of singular 
values, with respect to the background Lebesgue measure 
(hence with respect to more general invariant 
measures)?  If so, it might be interesting to develop a 
theory of random Hankel matrices in analogy with the 
theory of random Hermitian matrices.  

$(6.10)$ We have considered a number of invariant 
functions, $det(1+B(x)B(x)^{*})$, $det(1+W(\theta_{+})W(\theta_{+}
)^{*})$, 
$det(1+W(q_{+})W(q_{+})^{*})$, and so on.  It would be interesting 
to understand the Taylor series of these functions in a 
$G$-equivariant way.

\bigskip

\centerline{References}

\bigskip

[B] P. Billingsley, Convergence of Probability Measures, 
John Wiley and Sons (1968).

[Bo] V.I.  Bogachev, Gaussian Measures, A.M.S.  Math.  
Surveys and Monographs, Vol.  62 (1998).  

[FK] H. Farkas and I. Kra, Theta Constants, Riemann 
Surfaces and the Modular Group, Graduate Studies in 
Mathematics Vol. 37, AMS (2001).

[HH] W.  Helton and R.  Howe, Traces of commutators of 
integral operators, Acta Math.  135, No.  3-4 (1975) 
271-305.  

[He] S.  Helgason, Differential Geometry, Lie Groups, and 
Symmetric Spaces, Academic Press (1984).  

[Hi] E.  Hille, Analytic Function Theory, Vol.  II (1973).  

[GK], F.  Gardiner and N.  Lakic, Quasiconformal 
Teichmuller Theory, A.M.S.  Math.  Surveys and 
Monographs, Vol.  76 (2000).  

[K] J.P.  Kahane, Some Random Series of Functions, 
Cambridge Tracts in Advanced Mathematics (1985).  

[Kra] I.  Kra, Automorphic Forms and Kleinian Groups, 
Benjamin (1972).  

[L] O.  Lehto, Univalent Functions and Teichmuller 
Spaces, Springer Verlag (1986).  

[MKM] K. Matthes, J. Kerstan and J. Mecke, Infinitely 
Divisible Point Processes, John Wiley (1978).

[NS] S.  Nag and D.  Sullivan, Teichmuller theory and the 
universal period mapping via quantum calculus and the 
$H^{1/2}$ space on the circle, Osaka J.  Math.  32 (1995) 1-34.  

[Pe] V.  Peller, Hankel Operators and Their Applications, 
Springer-Verlag (2003).  

[Pi1] D.  Pickrell, Invariant measures for unitary groups 
associated to Kac-Moody Lie algebras, Memoirs of the 
A.M.S.,  No.  693 (2000).  

[Pi2] --------, An invariant measure for the loop space of 
a simply connected compact symmetric space, J.  Funct.  
Anal.  234 (2006) 321-363.  

[Pi3] --------, Loop spaces, diagonal distributions, and the 
Duistermaat-Heckman theorem, in progress.  

[Pi4] -------, A coordinate expression for the invariant 
measure of a loop group, in progress.
 
[PS] A.  Pressley and G.  Segal, Loop Groups, Oxford 
University Press (1986).  

[Se] I.E.  Segal, Ergodic subgroups of the orthogonal 
group on a real Hilbert space, Annals of Mathematics, 
Vol.  66, no.  2 (1957) 297-303.  

[St] D.  Stroock, Gaussian measures in traditional and not 
so traditional settings, Bull.  A.MS.,  Vol 33, No.  2 
(1996) 135-156.  

[T] W. Thurston, Three-dimensional Geometry and 
Topology, edited by S. Levy, Princeton Mathematical 
Series 35, Princeton University Press (1997).

[Z] R.  Zimmer, Ergodic Theory and Semisimple Groups, 
Birkhauser (1984).  

\end

There are four appendices; all are incomplete and merely 
attempt to frame interesting further questions.  

Appendix A discusses some background material on 
invariant (mainly normed) subspaces of $H^m$, and its 
formal completion.  

Appendix B concerns the Gaussian measure for 
$H^0/\Bbb C=H^1$, and the question of what $G$-ergodicity 
implies in terms of the dynamics of $G$.  The basic 
problem is to find a $G$-stable domain for the measure, 
where we can observe chaotic dynamics.  It is not clear 
that the abstract Wiener space approach is the right fit 
for this question.  For the particular example we are 
considering, this is closely related to the 
well-established subject of random Taylor series 
([Kahane]).  

Appendix C concerns the analysis of the Taylor expansion 
of $(0.12)$, in terms of the action of $G$.  It is not clear 
that this expansion has any immediate use.  However the 
combinatorics involved seem quite interesting, and this 
is an archetype for a large number of more complicated 
determinants.  

Appendix D concerns the $s$-number distribution of a 
random Hankel matrix.

\centerline{Appendix A. Invariant Subspaces}

\bigskip

The action of $\tilde {G}$ on $H^m$ is well-defined for any complex 
number $m$, but we will restrict attention to $m\in \Bbb R$.  We 
are mainly interested in the case $m>0$.  In this case 
the action is essentially irreducible, meaning that any 
invariant subspace is actually dense.  In this appendix 
we collect a number of elementary facts about these 
actions.  

The infinitesimal action of $sl(2,\Bbb C)$ on $H^m$ is given by 
$$A=\left(\matrix 0&0\\
1&0\endmatrix \right)\to -\frac {\partial}{\partial z}$$
$$A^{*}=\left(\matrix 0&1\\
0&0\endmatrix \right)\to z^2\frac {\partial}{\partial z}+2mz\tag A.1$$
$$H=\left(\matrix 1&0\\
0&-1\endmatrix \right)\to 2z\frac {\partial}{\partial z}+2m$$
The spectrum of the Hamiltonian is 
$$2m,2m+2,2m+4,..\tag A.2$$
For the infinitesimal action, the subspace $H^m_{alg}$, the span 
of the $z^n(dz)^m$ (the eigenfunctions of $H$), is an invariant 
subspace.  Also the infinitesimal action makes perfectly 
good sense on the dual, the formal completion $H^m_{formal}$ 
(where $f$ is allowed to be a formal power series).  

When $m$ is a nonpositive half-integer, the action is not 
irreducible.  We have already encountered two basic 
examples:  
$$0\to \Bbb C\to H^0@>{\partial}>>H^1\to 0\tag A.3$$
which is the linearization of the nonlinear mapping 
$$G\to H^0(\Delta ,G)\to H^1(\Delta ,G):g_{+}\to g_{+}^{-1}\partial 
g_{+},\tag A.4$$
and 
$$0\to sl(2,\Bbb C)\to H^{-1}@>{dS}>>H^2\to 0\tag A.5$$
which is the linearization of the Schwarzian 
$$SL(2,\Bbb C)\to exp(H^{-1})@>{S}>>H^2\tag A.6$$
where we think of $exp(H^{-1})$ as formal holomorphic 
automorphisms of the disk, whereas $SL(2,\Bbb C)$ is the 
complexification of the global holomorphic 
automorphisms). 

 In general, for $s=0,\frac 12,1,..$, a 
nonnegative half integer, for $m=-s$, 
$$0\to \Bbb C^{2s+1}=span\{z^0(dz)^{-s},..,z^{-2s}(dz)^{-s}\}\to 
H^{-s}\to H^{1-s}\to 0.\tag A.7$$
I do not know how to interpret this geometrically, 
except in the special cases above.  This might somehow 
explain why the cases $m=0$ or 1, and $m=2$ are so 
important.  

For each $m\ge 0$, the action $\tilde {G}\times H^m(\Delta )$ contains an 
irreducible invariant unitary action.  If $m>\frac 12$, then the 
essentially unique invariant Hilbert space norm is given 
by 
$$\vert f(dz)^m\vert^2_{H^m\cap L^2}=\frac {2m-1}{2\pi}\int_{\Delta}
\vert f(z)\vert^2(1-\vert z\vert^2)^{2m-2}dx\wedge dy=$$
$$\frac {2m-1}{2\pi}\sum_{n\ge 0}\vert f_n\vert^22\pi\int_0^1r^{2
n}(1-r^2)^{m-2}rdr=\sum_{n\ge 0}\vert f_n\vert^2\frac {\Gamma (n+
1)\Gamma (2m)}{\Gamma (2m+n)}$$
$$=\sum_{n\ge 0}\vert f_n\vert^2\frac {n!}{(2m)(2m+1)..(2m+n-1)}$$
$$=\vert f_0\vert^2+\frac 1{2m}\vert f_1\vert^2+\frac {2!}{(2m)(2
m+1)}\vert f_2\vert^2+..\tag A.8$$
where $f=\sum f_nz^n$.  The seemingly unnatural factor $\frac {2m
-1}{2\pi}$ 
is inserted because the last sum shows that this norm 
can be analytically continued to $m>0$, since the 
coefficients $B(n+1,2m)=\frac {\Gamma (n+1)\Gamma (2m)}{\Gamma (2
m+n)}\ge 0$.  In the critical 
case $m=\frac 12$, the norm can be understood as the $L^2$ norm 
of the $S^1$-boundary values of $f$ with respect to 
$\vert dz\vert =d\theta$.  For $0<m<\frac 12$, there is not an integral 
representation for the norm, but we will nonetheless 
continue to denote the norm by $\vert\cdot\vert_{H^m\cap L^2}$.  

The limit $m\downarrow 0$ is exceptional.  For $m=0$, for the norm 
defined as above, $\Bbb C$ (the constants) is an invariant 
subspace, with quotient determined by $\partial :H^0(\Delta )\to 
H^1(\Delta )$.  
However, in the definition of the norm, we could 
multiply by $m$.  In this case we obtain a Hilbert space 
substructure for the quotient $H^0/\Bbb C$ (the original vacuum 
is now a ghost), and an isometry 
$$H^0\cap L^2/\Bbb C@>{\partial}>>H^1\cap L^2:x\to\theta =\partial 
x.\tag A.9$$
We will leave the $m$ in the definition of the norm.  
However, especially in $\S 3$, we will make extensive use of 
the natural isomorphism $H^0/\Bbb C$ and $H^1$.  

We want to discuss the corresponding $L^p$ spaces.  We 
define $H^m\cap L^p$ to consist of $m$-differentials satisfying 
$$\vert f(dz)^m|_{L^p}^p=\int_{\Delta}\vert f(z)(dz)^m\vert_{\kappa}^
pdA_{hyperbolic}$$
$$=\int_{\Delta}\vert f\vert^p(1-\vert z\vert^2)^{mp-2}dx\wedge d
y<\infty .\tag A.10$$

In the case $p=\infty$ ...

\proclaim{(A.11)Proposition}This space is nontrivial if and 
only if $mp>1$.  
\endproclaim

\demo{Idea of Proof} 

In the critical case 
$$\int\vert f\vert^p(1-r^2)^{mp-2}d\theta rdr=\int_{u=0}^{\infty}
\int_{\theta}\vert f(\sqrt ue^{i\theta})\vert^pd\theta (1-u)^{-1}
du\tag A.12$$
The $\theta$ integral is an increasing function of $u$ [Rudin].  
\qed
\enddemo

Consider the critical case $mp=1$, or $p=\frac 1m$.  To obtain a 
nonempty space, we have to abandon the integral over $\Delta$, 
and replace it with an integral over the boundary, $S^1$.  
The resulting space $H^{1/p}\cap L^p$ is the classical Riesz 
space (which is classically denote by $H^p$, as in [Rudin], 
but of course we must avoid this notation).  The point 
is that 
$$\vert f(z)(dz)^m\vert^{1/m}=\vert f(z)\vert\vert dz\vert\tag A.13$$
can be integrated over $S^1$.  

In particular $H^0\cap L^{\infty}$ is the usual $H^{\infty}$ space of classical 
analysis (this is the only classical space which should 
naturally be thought of as consisting of functions!).  

We should ask the question of whether we can 
analytically continue this norm, as is possible in the 
case $p=2$.

\flushpar\ (A.14) Note:  The Poincare metric on $\Delta$:  
$$d(z,w)=2arctanh\vert\frac {z-w}{1-\bar {z}w}\vert .$$
So any function of this is Poincare invariant.  

\bigskip

\centerline{Appendix B. Analysis of The Taylor Expansion}

\bigskip

Consider the standard expansion 
$$det(1+B(x)B(x)^{*})=1+\vert B(x)\vert_{\Cal L_2}^2+\vert\Lambda^
2(B(x))\vert^2_{\Cal L_2}+..\tag B.1$$
as a function of $\theta\in H^1\cap L^2$.  Note that $tr(\Lambda^
n(BB^{*}))$ is 
homogeneous of degree $2n$ in $\theta$.  This means that the 
right hand side of $(B.1)$ is the Taylor expansion centered 
at $x=0$ for the left hand side as a function of $x$.  

We want to understand the meaning of these terms as 
$PSU(1,1)$-invariant functions.  For example we have 
already noted that in terms of $\theta$, the map 
$$\theta\to B(x)=\left(\matrix .&.\\
x_2&x_3&..\\
x_1&x_2&x_3&..\endmatrix \right)\tag B.2$$
can be understood PSU(1,1)-equivariantly as the linear 
map 
$$H^1\to H^{1/2}\otimes H^{1/2}=\Cal L_2(\bar {H}^{1/2}\otimes H^{
1/2})\tag B.3$$
where 

\proclaim{(B.4)Lemma}As a $SU(1,1)$ representation 
$$H^{1/2}\otimes H^{1/2}=H^1\oplus H^2\oplus H^3\oplus ..$$
More generally, 
$$H^s\otimes H^t=H^{s+t}\oplus H^{s+t+1}\oplus ..$$
Furthermore multiplication is the orthogonal projection 
$$0\to Ker\to H^s\otimes H^t@>{M}>>H^{s+t}\to 0.$$
and
$$S^2(H^{1/2})=H^1\oplus H^3\oplus ..$$
whereas 
$$\Lambda^2(H^{1/2})=H^2\oplus H^4\oplus ..$$
\endproclaim

\demo{Proof} The character of $H^s$ is $q^s/(1-q)$.  Now 
simply observe that the characters of the two sides are 
equal:  
$$\frac {q^s}{1-q}\cdot\frac {q^t}{1-q}=\frac {q^{s+t}+q^{s+t+1}+
..}{1-q}.\tag B.5$$
\qed
\enddemo

It is possible that $(B.1)$ above equals 
$$1+\vert\theta\vert^2+\vert P_2\theta^2\vert^2+\vert P_3\theta^3
\vert^2+..\tag B.6$$
where the projection in the quartic term corresponds to 
$$S^2(H^1)=H^2\oplus H^4\oplus H^6...\to\Lambda^2(H^{1/2})\otimes
\Lambda^2(H^{1/2})=$$
$$=(H^2\oplus H^4\oplus ..)^2=H^4\oplus ..\tag B.7$$
This injection apparently just drops the $H^2$ term.  Here 
we are using the amazing fact that 
$$S^2(H^1)=\Lambda^2(H^{1/2})\tag B.8$$
In general we are looking at (this is a first 
approximation - see below for the real decomposition of 
these spaces) 
$$S^n(H^1)=H^n\oplus H^{n+1}\oplus ..\to\Lambda^n(H^{1/2})\otimes
\Lambda^n(H^{1/2})$$
$$=(H^{n^2/2}\oplus H^{(n^2/2)+n}\oplus ..)\otimes (H^{n^2/2}\oplus 
H^{(n^2/2)+n}\oplus ..)$$
$$=H^{n^2}\oplus ..\tag B.9$$
where unfortunately some multiplicity begins to appear 
in the formulas.  But this does give us an idea of what 
we must cut off in the projections $P_n$, if this idea 
makes any sense at all.  

We now consider some general character calculations, 
which substantiate what we have been claiming.  

The character of $H^1$ is 
$$q+q^2+q^3+..\tag B.10$$
Therefore the character of $S(H^1)=S(q)\otimes S(q^2)\otimes ..$ equals 
$$\frac 1{1-q}\frac 1{1-q^2}...=(1+q+q^2+..)(1+q^2+q^{2\cdot 2}+.
.)(1+q^3+q^{3\cdot 2}+..)..$$
$$=1+q+2q^2+3q^3+5q^4+..$$
$$=\sum\vert \{(m_1,..,m_l):\sum m_jj=N\}\vert q^N=\sum_Np(N)q^N,\tag B.11$$
where the $m_i\ge 0$, and $p(N)$ is the partition function.  
The $m_i$ are the multiplicities for the corresponding 
partition:  the partition itself is given by 
$l^{(m_l)}>(l-1)^{(m_{l-1})}>..$, where the superscripts are 
indicating the number of occurrences of the elements in 
the partition.  

We now claim that the character of $S^n(H^1)$ is given by 
$$\sum\vert \{\vec {m}:\sum m_j=n,\sum m_jj=N\}q^N=\sum\vert \{\lambda
\in \Cal P(N):l(\lambda )=n\}\vert q^N;\tag B.12$$
here $\Cal P(N)$ denotes the set of partitions of $N$, and by 
definition the length of a partition is the number of 
nonzero terms in the partition, which is the same as the 
sum of the multiplicities.  For all of this see chapter 1 
of MacDonald.  To prove the claim, we go back to the 
calculation of the character, and the isomorphism of 
$S(H^1)$ with the tensor product of the $S(q^j)$.  To get 
something homogeneous of degree $n$ in $S(H^1)$, we start 
with elements homogeneous of degree $m_1$ in $S(q^1)$,..  and 
multiply them together.  The basic point is that the 
degree $m_i$ agrees with the multiplicity for the 
corresponding partition.  This proves the claim (in a 
rather inelegant way).  

We will denote the set of partitions of $N$ with length $n$ 
by $\Cal P_n(N)$, and its cardinality by $p_n(N)$.  

We can now write, as an isomorphism of group 
representations, 
$$S^n(H^1)=\sum (p_n(N)-p_n(N-1))H^N\tag B.13$$
Observe that there is an injective map 
$$\Cal P_n(N-1)\to \Cal P_n(N):(\mu_1\le ..\le\mu_n)\to (\mu_1\le 
..\le\mu_n+1).\tag B.14$$
The image of this map is exactly 
$$\{\lambda\in \Cal P_n(N):m(\lambda_n)=1\}\tag B.15$$
Therefore 
$$p_n(N)-p_n(N-1)=\vert \{\lambda\in \Cal P_n(N):m(\lambda_n)>1\}
\vert$$
$$=\vert \{\lambda\in \Cal P_n(N):\lambda_{n-1}=\lambda_n\}\vert\tag B.16$$
This is clearly going to be an arbitrarily large number, 
as $N\uparrow\infty$.  

\smallskip

\flushpar(B.17) Exceptional Example.  Suppose that $n=2$.  
Then 
$$\{\lambda\in \Cal P_2(N):\lambda_{n-1}=\lambda_n\}=\{\lambda =(\frac 
N2\le\frac N2)\}$$
Hence 
$$S^2(H^1)=H^2\oplus H^4\oplus H^6\oplus ..\tag B.18$$
If we consider $n=3$, then for $N=12$, there are at least 
two admissible partitions, $4\le 4\le 4$, and $2\le 5\le 5$.  Thus 
$S^3(H^1)$ is not multiplicity free.  
$$S^3(H^1)=H^3\oplus H^5\oplus H^6\oplus H^7\oplus H^8\oplus 2H^9
\oplus H^{10}\oplus 2H^{11}\oplus 2H^{12}\oplus ..\tag B.19$$

Here is a summary of what we have proven:  

\proclaim{(B.20)Proposition} The character of $S^n(H^1)$ 
equals 
$$\sum_Np_n(N)q^N.$$
As a representation of $PSU(1,1)$, $S^n(H^1)$ equals 
$$\sum_N\{p_n(N)-p_n(N-1)\}H^N$$
$$=\sum_N\vert \{\lambda\in \Cal P_n(N):m(\lambda_n)>1\}H^N$$
\endproclaim

The basic question:  how do we understand the natural 
map 
$$S^n(H^1)\to\Lambda^n(H^{1/2})\otimes\Lambda^n(H^{1/2})\tag B.21$$
Our first approximation to the truth:  look at the 
isotypic decomposition of $S^n$ and simply the first $n^2$-1 
possible isotypic components.  This is apparently true 
for $n=2$ because there is not any multiplicity, but 
beyond this we are in big trouble.  

\bigskip

Suppose that we consider $H^2$ in place of $H^1$.  

The character of $H^2$ is 
$$q^2+q^3+..\tag B.22$$
Therefore the character of $S(H^2)=S(q^2)\otimes S(q^3)\otimes ..$ 
equals 
$$\frac 1{1-q^2}\frac 1{1-q^3}...=(1+q^2+q^{2\cdot 2}+..)(1+q^3+q^{
3\cdot 2}+..)..$$
$$=1+q^2+q^3+2q^4+..$$
$$=\sum\vert \{(m_2,..,m_l):\sum m_jj=N\}\vert q^N,\tag B.23$$
and the character of $S^n(H^2)$ equals 
$$\sum\vert \{\lambda\in \Cal P(N):m_1(\lambda )=0,l(\lambda )=n\}
\vert q^N\tag B.24$$
and as a representation of $PSU(1,1)$, 
$$S^n(H^2)=\sum M(N)H^N,\tag B.25$$
where 
$$M(N)=\vert \{\lambda\in \Cal P_n(N):m_1(\lambda )=0\}\vert -\vert 
\{\lambda\in \Cal P_n(N):m_1(\lambda )=0\}\vert .\tag B.26$$

\bigskip

The character of $H^{1/2}$ equals 
$$q^{1/2}+q^{3/2}+..\tag B.27$$
Therefore the character of $\Lambda (H^{1/2})$ equals 
$$(1+q^{1/2})(1+q^{3/2})(1+q^{5/2})..\tag B.28$$
So the character of $\Lambda^2(H^{1/2})$ equals 
$$q^2+q^3+q^4+q^5+..$$
$$+q^4+q^5+....$$
$$+q^6+..\tag B.29$$
and as a representation 
$$\Lambda^2(H^{1/2})=H^2\oplus H^4\oplus H^6\oplus ..\tag B.30$$

\proclaim{(B.31)Proposition}As a representation of $SU(1,1)$, 
each of the representations $\Lambda^n(H^{1/2})$ is multiplicity 
free; more precisely 
$$\Lambda^3(H^{1/2})=H^{9/2}\oplus H^{9/2+3}\oplus H^{9/2+6}\oplus 
..$$
$$\Lambda^4(H^{1/2})=H^8\oplus H^{12}\oplus ..$$
and in general 
$$\Lambda^n(H^{1/2})=H^{n^2/2}\oplus H^{(n^2/2)+n}\oplus H^{(n^2/
2)+2n}\oplus ..$$

\endproclaim

\bigskip

\centerline{Appendix C. Ergodicity and Wiener Spaces.}

\bigskip

Given a real separable Hilbert space $H$, there is a 
Gaussian measure associated to this Hilbert space.  This 
measure is invariant in an appropriate sense with 
respect to the orthogonal group of $H$.  A basic result of 
Irving Segal states for that any group $G\subset O(H)$ with 
trivial finite dimensional unitary representations, the 
action of $G$ is ergodic ([Se]).  The point of this appendix 
is to understand what this means in terms of the 
dynamics of $G$, in what is probably the most elemental 
case, involving $G=PSU(1,1)$ acting on $H^1$.  

We first need to recall some background material 
concerning Gaussian measures.  An excellent starting 
point is [Stroock].  

If $P:H\to F$ is a finite dimensional orthogonal projection, 
and $E\subset F$ is a Borel subset, then $P^{-1}(E)$ is called a 
cylinder set.  The set of all cylinder sets, as $P$ varies, 
is an algebra.  

The first basic fact is that there is a unique finitely 
additive probability measure $\nu_H$ on the algebra of 
cylinder sets which is characterized by 
$$P_{*}\nu_H=\nu_F,$$
where $\nu_F$ denotes the standard Gaussian measure on $F$, 
i.e.  
$$\nu_F=(\frac 1{\sqrt {2\pi}})^{dim(F)}e^{-\frac 12\vert x\vert^
2}dm_F(x),$$
for all finite orthogonal projections $P$.  But it is 
impossible to extend $\nu_H$ to a countably additive measure 
on $H$.  

Suppose that $\vert\vert\cdot\vert\vert$ is a norm on $H$, and let $
B$ denote the 
Banach space completion.  A cylinder set for $B$ is of the 
form $P^{-1}(E)$, where $P$ is an orthogonal projection as 
above, but we now additional require that $P$ has a 
continuous extension to $B$, and we regard $P:B\to F$.  In the 
same way as we constructed $\nu_H$, we can construct a 
finitely additive probability measure $\nu_B$.  

\proclaim{Summary of Basic Facts} 

(a) [Gross, Dudley-Feldman-LeCam] $\nu_B$ extends to a 
countably additive measure on $B$ if and only if the norm 
is measureable in the sense of Gross:  .  

Assume that $\nu_B$ does have a countably additive extension.  

(b) [Cameron-Martin, Segal] $\nu_B$ is $H$-translation 
quasiinvariant, and positive on all open subsets of $B$; in 
particular it is impossible for $\nu_B$ to be supported on a 
sphere of $B$.  

(c) [Segal] $\nu_B$ is orthogonally invariant in the following 
senses:  

(1) If we consider the subgroup of $O(H)$ consisting of 
transformations which extend continuously to maps $B\to B$, 
then $\nu_B$ is invariant under this subset.  

(2) The corresponding unitary representation has a 
unique continuous extension 
$$O(H)\to U(L^2(\nu_B))$$
where $O(H)$ has the norm topology.  

(d) Given a locally compact second countable $G\subset O(H)$, it 
is possible to choose $B$ so that $G$ acts continuously on $B$ 
by bounded operators.  

(e) If $G$ as in (d) has the property that all of its finite 
dimensional unitary representations are trivial, then $\nu_B$ 
is $G$-ergodic.  As a consequence of $(b)$ and ergodicity, for 
a.e.  $b\in B$, the $G$-orbit $G\cdot b$ is dense in $B$.  

(f) It is impossible to choose $B$ so that $G$ acts 
isometrically.  
\endproclaim

I do not have a proof of $(d)$.  

We would like to follow up on $(e)$.  How can we describe 
the $b$ so that the orbit is dense.  This is where we have 
to understand what a typical $b$ is.  

Let's now consider the particular case of $G=PSU(1,1)$ 
acting on $H^0/\Bbb C\cap L^2=H^1$ and the Gaussian $\nu_T=\nu^{(
1)}_T$.  Can 
we 

(1) find a really nice $B$ that is $G$-stable, but is really 
tight.  

(2) Can we use the results from [Kahane] to describe $b$ 
with dense orbits.  

\smallskip

Remark.  Useful Fact:  Suppose that $\mu$ is a standard 
measure, positive on open sets.  Given an ergodic action 
$G\times (x,\mu )$, for a.e.  x, $G\cdot x$ is dense.  Zimmer, page 10.  

\smallskip

We want to define a conformally invariant norm of the 
tentative form 
$$\vert x\vert_B=\sup\frac {\vert x(z)-x(w)\vert}{\sqrt {d(z,w)ln
(d(z,w))}}$$
Note that the values involved in this $\sup$ are small 
when $z$ and $w$ are close, so the larger values tend to 
occur when $z$ and $w$ are separating.  This is motivated 
by Levy's Law, for Brownian motion, and we are trying 
to mimic this in this conformally invariant context.  

We would like to prove something like 
$$\limsup_{d(z,w)\uparrow\infty}\frac {\vert x(z)-x(w)\vert}{\sqrt {
dln(d)}}=\sqrt T$$
whereas we know something like 
$$\limsup_{r\uparrow 1}\frac {\vert x(re^{i\alpha})\vert}{\sqrt {
ln(\frac 1{1-r^2})ln(\frac 1{1-r}ln(\frac 1{1-r}))}}=\sqrt T$$

In connection with Kahane's formulas in \S 2, recall that 
$$arctanh(r)=\frac 12(ln(1+r)+ln(\frac 1{1-r}))$$
and 
$$d(z,w)=2arctanh\frac {\vert z-w\vert}{\vert 1-\bar {z}w\vert}.$$

\bigskip

\centerline{Appendix D. Isospectral Sets and $s$-number distributions}

\bigskip

Given a Hilbert-Schmidt Hankel operator $B(x)$, let 
$s_1\ge s_2\ge ..$ denote the singular values of $B$, i.e.  the 
eigenvalues of $\sqrt {BB^{*}}$.  We think of $s_j$ as a function of 
$\theta\in H^1\cap L^2$, or $x\in H^0/\Bbb C\cap L^2$.  If we have a conformally 
invariant measure on $H^1$, it is natural to ask about the 
distribution of the $s$-numbers of $B_n$ as $n\to\infty$.  We have 
also had to deal with integrals of functions which depend 
only on the $s$-numbers, so again we need to ask about 
the distribution.  

\smallskip

\flushpar Example.  Consider the case $n=2$.  In terms of 
singular values, the integral $()$ of \S 3 can be written as 
$$\int\frac 1{(1+s_1^2)^p(1+s_2^2)^p}dm(x_1,x_2)$$
The eigenvalues of $B_2B_2^{*}$ are 
$$s_1^2,s_2^2=\frac 12(tr\pm ((tr)^2-4det)^{1/2}=$$
$$\frac 12\{(\vert x_1\vert^2+2\vert x_2\vert^2)\pm ((\vert x_1\vert^
2+2\vert x_2\vert^2)^2-4\vert x_2\vert^4)^{1/2}\}$$
Let $u_i=\vert x_i\vert^2$, $\sigma_2=det(BB^{*})$, $\sigma_1=tr(
BB^{*})$, $S_i=s_i^2$.  
Because 
$$\sigma_1=u_1+2u_2=S_1+S_2,$$
$$\sigma_2=u_2^2=S_1S_2,$$
$$du_1\wedge du_2=d\sigma_1\wedge d(\sigma_2^{1/2})=\frac {d\sigma_
1\wedge d\sigma_2}{2\sigma_2^{1/2}}=\frac {(S_2-S_1)dS_1\wedge dS_
2}{2(S_2S_1)^{1/2}}$$
$$=(s_2^2-s_1^2)ds_1\wedge ds_2$$

The example $()$ is then a Selbert integral (chapter of 
[Mehta]), and reproduces $()$.  

A basic general problem is to compute the distribution of 
the $s_j$, as we have done here.  

In general it is a basic fact that 
$$d\sigma_1\wedge ..\wedge d\sigma_n=\prod_{i<j}(S_j-S_i)dS_1\wedge 
..\wedge dS_n$$
where $S_i=s_i^2$ and the $\sigma_i$ are the elementary symmetric 
functions in the $S_i$.  So we can just as well consider 
the distribution of the $\sigma_i$.  

Consider the case $n=3$.  The symmetric functions of the 
singular values are 
$$\sigma_1=\vert x_1\vert^2+2\vert x_2\vert^2+3\vert x_3\vert^2,\quad$$
$$\sigma_2=\vert x_1x_3-x_2^2\vert^2+2\vert x_2x_3\vert^2+3\vert 
x_3\vert^4$$
$$=\vert x_1x_3\vert^2+\vert x_2\vert^4-2\Bbb Re(x_1x_3\bar {x}_2^
2)+2\vert x_2x_3\vert^2+3\vert x_3\vert^4,$$
$$\sigma_3=\vert x_3\vert^6.$$
Let $X_i=\vert x_i\vert^2$.  
$$d\sigma_1=dX_1+2dX_2+3dX_3,$$
$$d\sigma_2=X_1dX_3+X_3dX_1+2X_2dX_2-($$
$$2d\Bbb Re(x_1x_3\bar {x}_2^2)=[X_2(X_1X_3)^{-1/2}dX_1+2(X_1X_3)^{
1/2}dX_2]cos(..))+$$
$$2(X_2dX_3+X_3dX_2)+6X_3dX_3+$$
stuff that involves $dX_3$ and $d\theta$ 
$$d\sigma_3=3X_3^2dX_3$$
Therefore 
$$d\sigma_1\wedge d\sigma_2\wedge d\sigma_3\wedge d\theta =$$
$$=(dX_1+2dX_2)\wedge ((X_3-\frac {X_2}{\sqrt {X_1X_3}}cos)dX_1+(
2X_2-2\sqrt {X_1X_3}cos)dX_2)$$
$$\wedge 3X_3^2dX_3\wedge d\theta$$
$$=2(X_2-\sqrt {X_2X_3}cos-X_3+\frac {X_2}{\sqrt {X_1X_3}}cos)X_3^
2dX\wedge d\theta$$
We need to express this density in terms of $\sigma$ and $\theta$, 
which seems impossible.  

We have invariant functions 
$$tr(\Lambda^n(B)\Lambda^n(B)^{*}),\quad n=1,2,3,..$$
In this particular case these invariant functions are 
actually the elementary symmetric functions in the 
$s$-numbers of $B$:  
$$det(1+BB^{*})=1+tr(BB^{*})+..$$
$$=\prod (1+s_j(B)^2)=1+\sum s_j^2+..$$
$$=\sigma_0+\sigma_1(s_1^2,..)+\sigma_2(s_1^2,..)+..$$
When we are considering infinitely many variables, it is 
often useful to also consider the Schur functions, which 
form an orthonormal basis for the $\Bbb Z$-ring of universal 
symmetric polynomials (see page 202 of [PS] and 
[MacDonald]).

$======================$

\bigskip

\centerline{Section 3.}

\smallskip

Can we obtain the Hankel measures from transitive 
measures?  For example suppose that we can 
holomorphic one forms $\theta_j$ corresponding to surfaces of 
genus $j$.  We could say something more precise about 
the choice of the $\theta_j$ , such as that they correspond to 
an $\alpha$ cycle, or something more complicated.  We then 
have transitive measure $\nu_1$, $\nu_2$, ..  What can we say 
about the convolution of these guys?  

We could ask a similar question for $H^2$ for more 
generally $H^m$ when $m$ is integral.

\centerline{Universality Questions}

\smallskip

Suppose that we fix $k$ and $n$.  Consider the question of 
whether 
$$\int\frac 1{(1+\vert\Lambda^k(B_n)\vert^2)^p}dm(x_1,..,x_n)<\infty$$
Is it possible that there is a limiting measure for each 
$k$, together with a deformation parameter $l$?  

Suppose that $k>1$.  If $n=1$, then $\Lambda^k(B_1)=0$.  Hence 
$p_1=\infty$.  This strongly suggests that this fails.  

So we would need to consider something like 
$$\int\frac 1{(\sum_0^k\vert\Lambda^j(B_n)\vert^2)^p}dm(x_1,..,x_
n)$$
Now we will have critical exponents $p_n(k)$.  

In the case $k=1$, our favorite, $p_n(1)=n$.  So our 
measure looks like 
$$\frac 1{\Cal Z_n}\frac 1{(1+\vert x\vert^2)^{n+1+l}}dm(x_1,..,x_
n)$$
where $-1<l$.  Very familiar.  

The question is what can we say about critical 
exponents for $1<k<\infty$.  We could look at a table, and 
there might be some interesting pattern.  I think we 
should expect that for fixed $k$, asymptotically $p_n(k)$ 
behaves as in the case $k=1$, $p_n(k)\tilde{}$$n$.  What we know so 
far is summarized as follows:  

$$\matrix k=&1&2&3&..&\infty\\
n=1&1&1&1&&2-1\\
n=2&2&2-\frac 12&2-\frac 12&..&2-\frac 12\\
n=3&3&?&2-\frac 13&..&2-\frac 13\\
n=4&4&?&?&&2-\frac 14\\
n=5\endmatrix $$

Now suppose $1<k$.  The next question is whether 
$$\lim_{n\uparrow\infty}\frac 1{\Cal Z_n}\frac 1{(\sum_0^k\vert\Lambda^
j(B_n)\vert^2)^{p_n(k)+L}}dm(x_1,..,x_n)$$
exists, and if it exists, is it just a linear measure, or is 
it something new?  Our conjecture above for the critical 
exponents anticipates that, despite appearances, this 
would tend to the linear case.  So there would be a 
dichotomy.  

Something we are totally missing at this point is a 
measure which has $x_1$-density 
$$\frac 1{(1+\vert x_1\vert^2)^{3/2}}F(\frac 1{1+\vert x_1\vert^2}
)$$
where $F(\rho )$ is the mysterious function 
$$F(\rho )=\int_0^{\infty}\frac {\rho}{(\rho +\frac {(x-1)^2}x)^{
3/2}}dx$$
The $F(\rho )$ does not affect the qualitative properties of 
the measure.  It is entirely possible that we have the 
two measures mixed up.  

Note that there is a very simple Fourier transform 
$$(\frac 1{(1+\vert x_1\vert^2)^{3/2}})^{\hat{}}(y_1)=e^{-\vert y_
1\vert}$$
ignoring some constants here and there (Note we are 
dealing with a superposition of Gaussians; see 26.  page 
262 in [Folland]).  Convolving with $F^{\hat{}}$ will not change 
things too much.  

But these two measures, whatever they are, should be 
generators in some sense, for an interesting class of 
measures.  

\bigskip
\centerline{Section 4.}

\smallskip

Is $\frac 1{det(1+WW^{*})}$ a positive definite function on 
$H^1\cap L^2\otimes \frak k^{\Bbb C}$?  If we think about the $SU
(2)$ case, note that 
the restriction{\bf s} to $\alpha${\bf ,} $\beta${\bf ,} and $\gamma$ are positive definite.  But 
this does not automatically mean the whole thing is.  To 
find the corresponding measure one is tempted to look at 
the bracket 
$$H^{1/2}\otimes \frak k^{\Bbb C}\otimes H^{1/2}\otimes \frak k^{
\Bbb C}\to H^1\otimes \frak k^{\Bbb C}$$
and the image of the Gaussians.  

It would seem more natural to consider curvature, but 
this lands in $H^2$.  

\bigskip

\centerline{Section 5.}

\smallskip

In this case there is not a nice relationship between the 
$L^2$ properties of $q\in H^2$ and $W(q)$.  

We need to spell out the form of $W(q)$.  

Is there some connection between Loewner process and 
what we are discussing in this paper?  There might be 
some connection with the boundary of Teichmuller space.  

We are supposed to have $c=\frac {(3\kappa -8)(6-\kappa )}{2\kappa}$, $
h=\frac {6-\kappa}{2\kappa}$ The 
trace is a simple curve for $\kappa\le 4$, it has simple double 
points up to $8$, and it is space filling after that.  

For the radial Loewner process, we can take $\lambda_tu_t=g_t^{-1}$, 
$\lambda_t=e^{-t}$, $l=z$ on $\Delta^{*}$.  

\bigskip

For transitive measures, it is possible that it is never 
the case that their Fourier transforms have this 
continuity property.  This question, for transitive 
measures, might be equivalent to $\theta$ belonging to $L^2\cap H^
m$ 
(for a measure with finite dimensional support, this 
might be equivalent to literally having the support in 
$L^2\cap H^m$).  At least for $m\ge 1/2$, this condition is that 
$$\int_{\Delta}\vert\theta\vert^2dA(q)<\infty .$$
However the integrand is invariant under the lattice $L$.  
So this integral is infinite, because $L$ is infinite.  

The Fourier transform is given by 
$$\int_{L\backslash G}e^{-i\Bbb Re\langle\eta ,g^{*}\theta\rangle_
m}dV(g)$$
What we argued above is that no meaning can be 
attached to $\langle\eta ,g^{*}\theta\rangle$.  The correct way to think about this 
is to first think about this as a function which is 
well-defined on $(H^m)^{*}$.  Then we ask whether there 
exists a continuous extension to the larger space, the 
dual of $H^m\cap L^2$.  

I am uncertain about this argument, because if one 
considers the Gaussian, this is not supported on $H^m\cap L^2$, 
and yet its Fourier transform does have a continuous 
extension to the dual of this Hilbert space.  

This might be a more precise way of describing the 
class of measures that we are interested in.  

Our measures have the property that their Fourier 
transforms are continuous on the dual of $H^m$.  The 
measures that we are really interested in have Fourier 
transforms which have continuous extensions to the dual 
of $H^m\cap L^2$ (and which are invariant).  

Suppose that $\nu$ is a transitive measure.  Let $\theta_0$ denote a 
basepoint in the support, let $L$ denote the projection to 
$G$ of the stability subgroup, and let $\Cal F$ denote a 
fundamental domain for the action of $L$ on $\Delta$.  We can 
write the Fourier transform as 
$$\hat{\nu }(\theta )=\int e^{-i\Bbb Re\int_{\Delta}g^{*}\theta\wedge 
*\theta_0}d\nu$$
The inner product does not make any sense for $\theta$ a 
general $L^2$ form, because the integral is invariant under 
the lattice.  This strongly suggests the Fourier 
transform is not continuous on $H^m\cap L^2$ (but this requires 
more attention).

Can we possibly formulate some kind of multiplicative 
central limit theorem.  Consider a nonvanishing 
differential, such as $\eta^2dz$.  Let 
$$\theta_n=\eta^{2/n}(dz)^{1/n},\quad\nu_n=\delta_{G_{*}\theta_n}$$
We could ask about 
$$\lim_{n\to\infty}\nu_n\otimes ..\nu_n\in Prob(H^1).$$

An alternative is to consider dividing Gaussians, 
following the idea of Gelfand et al of continuous product 
of representations (the difference here is that Gelfand et 
al consider the complementary series, which does truly 
converge to the trivial representation.  If we consider 
$H^m$, this converges to $\Bbb C$ and $H^1$ as $m\downarrow 0$ (look at the 
character).  Given $m_1+..+m_n=1$ and $t_1,..,t_n$ we could 
consider 
$$\nu^{(m_1)}_{t_1}\otimes ..\otimes\nu^{(m_n)}_{t_n}$$
a probability measure on $H^1$.  Is there some natural way 
to obtain a limit from this.  We need to have some idea 
of how to think about the Fourier transform of this 
measure.  
 
Note that 
$$\delta_0\otimes\mu =\delta_0$$
So if the $t_i$ are getting small, we might just be going to 
zero.

\bigskip

About Fourier transforms.  

\bigskip

c) It is important here that the rank is one.  It is not 
the case that the invariant measure for a higher rank 
Grassmannian is a convex combination of Gaussians.  We 
have to look at the more complicated classification of 
ergodic invariant measures for the action of the 
stability subgroup acting on $Z\in \Cal L(H^{+},H^{-})$.  I wonder if 
Olshanksy and Borodin have actually worked out what 
the decomposition is in the infinite rank case.  What is 
it in the finite rank case?  This example suggests that 
the Gaussians cannot be ignored!

The first remark involves the standard isomorphism 
$$\hat {\Bbb C}\to \Bbb P(\Bbb C^2):z\in \Bbb C\to \Bbb P(\left(\matrix 
z\\
1\endmatrix \right)),\infty\to \Bbb P(\left(\matrix 1\\
0\endmatrix \right)).$$
The linear fractional action $()$ of $G$ on $\hat {\Bbb C}$ is a coordinate 
expression for the action of $G$ on $\Bbb P(\Bbb C^2)$, which is 
induced by the defining representation of $G^{(2)}=SU(1,1)$ 
on $\Bbb C^2$.  This displays the importance of the double 
cover.  

The second remark is that the holomorphic tangent space 
at a point $L\in \Bbb P(\Bbb C^2)$ is naturally identified with the 
complex line 
$$Hom(L,\Bbb C^2/L);$$
here natural can be interpreted to mean that the action 
of $PSL(2,\Bbb C)$ on the holomorphic tangent bundle is the 
obvious action on 
$$\bigsqcup_LHom(L,\Bbb C^2/L).$$
The canonical bundle $\kappa$, which is dual to the holomorphic 
tangent bundle, has fiber at $L$ which is naturally 
identified with the complex line 
$$(L^{*}\otimes (\Bbb C^2/L))^{*}=L\otimes L.$$
In terms of the isomorphism $()$, and the coordinate $z$, 
$$(dz)\vert_{z_0}\to$$

Now consider the $G$-stable decomposition $()$.  to think 
about these bundles is as follows.  Consider the canonical 
bundle $\kappa$ of $\hat {\Bbb C}$.  Up to isomorphism there is a unique 
holomorphic line bundle $\kappa^{1/2}$ such that its square is $\kappa$.  
In terms of the isomorphism $\kappa^{1/2}$ can be identified with 
the tautological holomorphic line bundle.  We have 
$$H^0(\kappa^{1/2}\vert_{\Delta})=H^{1/2}\quad\text{\rm and}\quad 
H^0(\kappa^{1/2}\vert_{\Delta^{*}})=H^{1/2}(\Delta^{*}).$$

The restriction of $\kappa^{1/2}$ to $S^1$ has a real structure.  This 
underlying real bundle is a realization of the odd spin 
structure.  The sections are naturally identified with 
odd functions on the circle, using the trivialization 
$(\frac {dz}z)^{1/2}$.  Let $\Cal D^{(2)}$ denote the double cover of $
\Cal D$, the 
group of orientation-preserving diffeomorphisms.  There 
is a natural action of this double cover on $\kappa^{1/2}\vert_{S^
1}$, and 
there is a $\Cal D^{(2)}$-invariant Hermitian inner product on the 
space of sections (the space of odd spinors), 
$$\Omega^0(\kappa^{1/2})\otimes\widebar{\Omega^0(\kappa^{1/2})}\to 
\Bbb C:\theta\otimes\bar{\eta}\to\int_{S^1}\theta\bar{\eta}$$
The point is that if $\theta =f(z)(dz)^{1/2}$ and $\eta =g(z)(dz)^{
1/2}$ 
are odd spinors on $S^1$, then $\theta\bar{\eta }=f\bar {g}\vert 
dz\vert$ is naturally a 
one density on $S^1$, which can be integrated.

\demo{Sketch of Proof} The careful proof of this will 
appear in another paper.  Let 
$$B_n=B(\sum^nx_jz^j).$$
and $\delta_n=det(1+B_nB_n^{*})$.  For purposes of calculating 
determinants, we can view $B_n$ as an $n\times n$ matrix (with 
$x_n$ along the main diagonal).  As an $(n+1)\times (n+1)$ matrix 
equation, we can write 
$$B_{n+1}=x_{n+1}1_{n+1}+b_n,$$
where $1_{n+1}$ is the $n+1\times n+1$ `identity matrix' and $b_n$ 
denotes $B_n$ viewed as a strictly lower triangular 
$n+1$$\times n+1$ matrix.  Then 
$$\delta_{n+1}=det(1_{n+1}+B_{n+1}B_{n+1}^{*})=$$
$$=det((1+\vert x_{n+1}\vert^2)1_{n+1}+x_{n+1}B_n^{*}+B_nx_{n+1}^{
*}+B_nB_n^{*})$$

To understand the main source of difficulty, suppose 
that we could simply ignore the interaction terms 
involving $x_{n+1}$ and $B_n$.  Then 
$$\int\frac 1{det(1+B_{n+1}B_{n+1}^{*})^p}d\lambda (x_{n+1})=$$
$$\int\frac 1{det(1+x_{n+1}b_n^{*}+b_nx_{n+1}^{*}+b_nb_n^{*}+\vert 
x_{n+1}\vert^2)^p}d\lambda (x_{n+1})$$
$$=\int_0^{\infty}\frac 1{\prod_{j=1}^n(1+u+s_j(B_n)^2)^p}\frac 1{
(1+u)^p}du$$
$$\le\frac 1{det(1+B_nB_n^{*})^p}\int_0^{\infty}\frac 1{(1+u)^p}d
u$$
The integral is finite provided that $1<p$.  If we could 
ignore the interaction, we could use an induction 
argument to claim that our critical $p=1$, which is not 
correct.  

What we have done, in ignoring the interaction, is kill 
all offdiagonal terms, as we will see through examples.  

$===========================$ 

To get the correct behavior, especially at infinity, we 
have to replace the exponent $p=2$ by $p=5/2$; we have 
to reference everything with respect to the critical 
exponent.  If we want to obtain the $x_1$ distribution with 
$p=2=p_1+1$, then we need to take the $x_1,x_2$ distribution 
with $p=p_2+1=5/2$.  

Try again:  
$$\int_{x_2\in \Bbb C}\frac 1{(\vert x_1\vert^2+(1+\vert x_2\vert^
2)^2)^{5/2}}dm(x_2)$$
$$=\int_{v=1}^{\infty}\frac 1{(\vert x_1\vert^2+v^2)^{5/2}}dv=\frac {
\vert x_1\vert}{\vert x_1\vert^5}\int_{u=\vert x_1\vert^{-1}}^{\infty}\frac 
1{(1+u^2)^{5/2}}du$$
$$=\frac 1{\vert x_1\vert^4}(\frac {3u+2u^3}{3(1+u^2)^{3/2}})\vert_{
u=\vert x_1\vert^{-1}}^{\infty}$$
$$=\frac 1{3\vert x_1\vert^4}(2-\frac {3\vert x_1\vert^{-1}+2\vert 
x_1\vert^{-3}}{(1+\vert x_1\vert^{-2})^{3/2}}).$$
This now has the correct behavior at infinity.  
Furthermore, if we normalize, the Taylor series is close 
to that for $1/(1+r^2)^2$:  
$$1-\frac 53r^2+\frac {35}{16}r^4-\frac {21}8r^6+\frac {385}{128}
r^8+..$$
versus 
$$1-2r^2+3r^4-4r^6+..$$
where $r=\vert x_1\vert$.  We are really quite close.  

Suppose that we consider the critical exponent.  
$$\int\frac 1{(\vert x_1\vert^2+(1+\vert x_2\vert^2)^2)^{3/2}}dm(
x_2)$$
$$=\int_{v=1}^{\infty}\frac 1{(\vert x_1\vert^2+v^2)^{3/2}}dv=\frac {
\vert x_1\vert}{\vert x_1\vert^3}\int_{u=\vert x_1\vert^{-1}}^{\infty}\frac 
1{(1+u^2)^{3/2}}du$$
$$=\frac 1{\vert x_1\vert^2}\frac u{(1+u^2)^{1/2}}\vert_{u=\vert 
x_1\vert^{-1}}^{\infty}$$
$$=\frac 1{\vert x_1\vert^2}(1-\frac {\vert x_1\vert^{-1}}{(1+\vert 
x_1\vert^{-2})^{1/2}}).$$
Sure enough, this does resemble $1/(1+r^2)$.  

Suppose $n=3$.  Write $x_3=u^{1/2}exp(i\theta )$.  Then 
$$\int\frac 1{\delta_3^p}dm(x_3)=$$
$$\int\frac {dm(x_3)}{(1+\vert x_1\vert^2+2\vert x_2\vert^2+3\vert 
x_3\vert^2+2\vert x_2x_3\vert^2+\vert x_3x_1-x_2^2\vert^2+3\vert 
x_3\vert^4+\vert x_3\vert^6)^p}$$
$$=\int \{\int\frac 1{(a-bcos(\theta ))^p}d\theta \}du$$
where 
$$a=1+\vert x_1\vert^2+2\vert x_2\vert^2+3u+2\vert x_2\vert^2u+\vert 
x_1\vert^2u+\vert x_2\vert^4+3u^2+u^3$$
$$=(\vert x_1\vert^2+2\vert x_2\vert^2)(1+u)+\vert x_2\vert^4+(1+
u)^3$$
$$b=2u^{1/2}\vert x_1\vert\vert x_2\vert^2$$
If $p=2,3,4,..$ we know how to evaluate this $\theta$ integral.  
We consider the case $p=3$:  using 
$$\int_0^{\pi}\frac 1{(1-Acos(t))^3}dt=\frac {\pi (2+A^2)}{2(1-A^
2)^{5/2}},\quad A^2<1$$
we see that the $x_1,x_2$ distribution equals 
$$\int_0^{\infty}\frac {2a^2+b^2}{(a^2-b^2)^{5/2}}du$$
Unfortunately Maple cannot cope with this.  

We should try $p=2$.  In this case the $x_1,x_2$ distribution 
is given by 
$$\int\frac 1{\delta_3^2}dm(x_3)=\int_0^{\infty}\frac a{(a^2-b^2)^{
3/2}}du$$
$$=\int_1^{\infty}\frac {(\vert x_1\vert^2+2\vert x_2\vert^2)v+\vert 
x_2\vert^4+v^3}{(((\vert x_1\vert^2+2\vert x_2\vert^2)v+\vert x_2
\vert^4+v^3)^2-4(v-1)\vert x_1\vert^2\vert x_2\vert^4)^{3/2}}dv$$
(exactly as in the $SU(3)$ case in the first appendix, 
except that $a$ is now a cubic in $u$, whereas previously it 
was quadratic).  

$==========================$ 

Let's temporarily give up on this and try to integrate 
the whole thing, beginning with $\vert x_1\vert$:  
$$\int\frac 1{((\vert x_1\vert^2+2\vert x_2\vert^2)(1+\vert x_3\vert^
2)-2\vert x_3\vert\vert x_1\vert cos(\theta )+\vert x_2\vert^4+(1
+\vert x_3\vert^2)^3)^p}$$
$$=\int\frac 1{((1+\vert x_3\vert^2)(\vert x_1\vert -\frac {\vert 
x_3\vert cos(\theta )}{1+\vert x_3\vert^2})^2-\frac {\vert x_3\vert^
2cos^2(\theta )}{1+\vert x_3\vert^2}+2\vert x_2\vert^2(1+\vert x_
3\vert^2)+\vert x_2\vert^4+(1+\vert x_3\vert^2)^3)^p}$$
Using 
$$\int_0^{\infty}\frac {rdr}{(\alpha (r-\beta )^2+\gamma )^p}=\frac 
1{\alpha (1-p)}\frac 1{(\alpha (r-\beta )^2+\gamma )^{p-1}}\vert_{
r=0}^{\infty}$$
$$=\frac 1{(p-1)}\frac 1{\alpha (\alpha\beta^2+\gamma )^{p-1}},$$
this is wrong - you forgot the $\beta$ involved in the change 
of variable.  

we see $*$ is equivalent to 
$$=\int\frac 1{(1+\vert x_3\vert^2)(2\vert x_2\vert^2(1+\vert x_3
\vert^2)+\vert x_2\vert^4+(1+\vert x_3\vert^2)^3)^{p-1}}$$
$$=\int\frac 1{(1+\vert x_3\vert^2)^{2-p}((1+\vert x_3\vert^2)(2\vert 
x_2\vert^2+\vert x_2\vert^4)+(1+\vert x_3\vert^2)^4)^{p-1}}$$
It seems hopeless to do this exactly for general $p$.  
However for $p=2$, we see this equals (with $v=\vert x_2\vert^2$) 
$$=\int\frac 1{(((1+\vert x_3\vert^2)^2-1)(\vert x_3\vert cos(\theta 
))^2+(1+\vert x_3\vert^2)(2\vert x_2\vert^2+\vert x_2\vert^4)+(1+
\vert x_3\vert^2)^4)^{p-1}}$$
$$=\frac {arctan(\frac {2b+2bv}{2\sqrt {ab}})}{\sqrt {ab}}\vert_{
v=0}^{\infty}$$
$$=\frac 1{(((1+\vert x_3\vert^2)^2-1)(\vert x_3\vert cos(\theta 
))^2-(1+\vert x_3\vert^2)+(1+\vert x_3\vert^2)^4)}$$
It seems clear that this will be finite when we integrate 
wrt $\vert x_3\vert^2$ and $\theta$.  So we definitely have $p_3\le 
2$, and it 
seems very plausible that $p_3<2$, but we probably cannot 
do this exactly.

We suppose that $\int\delta_n^{-1-\epsilon}<\infty$.  

Now 
$$\delta_{n+1}=\delta_n+(n+1)\vert x_{n+1}\vert^2+2\sum^n_1j\vert 
x_j\vert^2\vert x_{n+1}\vert^2+..+\vert x_{n+1}\vert^{2(n+1)},$$
where the intermediate terms are positive.  Hence 
$$\int\frac 1{\delta_{n+1}^{1+\epsilon}}dm(x_{n+1})=\int\frac 1{(
\delta_n+..+\vert x_{n+1}\vert^{2(n+1)})^{1+\epsilon}}dm(x_{n+1})$$
$$=\frac 1{\delta_n^{1+\epsilon}}\delta_n^{\frac 1{n+1}}\int\frac 
1{(1+...+\vert x_{n+1}\delta_n^{-\frac 1{2(n+1)}}\vert^{2(n+1)})^{
1+\epsilon}}dm(x_{n+1}\delta_n^{-\frac 1{2(n+1)}})$$
The intermediate terms which we have surpressed are 
of the form 
$$(n+1)\vert x_{n+1}\delta_n^{-\frac 12}\vert^2+..$$

\qed
\enddemo

There is an action of $\Bbb T^3$ which acts on the isosingular 
sets:  
$$(x_i)\to\left(\matrix \lambda_1\lambda_2\lambda_3^2x_1\\
\lambda_1\lambda_3x_2\\
\lambda_1\lambda_2^{-1}x_3\endmatrix \right)$$
This is not a faithful action, because $(\lambda ,\lambda ,\lambda^{
-1})$ acts 
trivially.  So we only have a $\Bbb T^2$ action.  This is just 
the $Rot(S^1)\times T$ action.  

It seems highly unlikely that there is any kind of 
integrability for $n>2$.  

If we assume that the $x_j$ are real, then it is easy to 
actually solve for the $x_j$ in terms of the $\sigma_j$.  

We have $x_3=\sigma_3^{1/6}$.  We then have a quadratic equation 
for $x_2^2$.  We then have a linear equation for $x_1$.  

If $n=4$ we may not be so lucky.  We can find $x_4$.  But 
then we have a cubic equation for $x_3^2$.  

But in this real case the simple thing to do is just to 
directly compute 
$$d\sigma_1\wedge d\sigma_2\wedge ..d\sigma_n=\delta dx_1\wedge .
.\wedge dx_n$$
This would be completely straight forward.  We could 
then solve for $dx$.  

In the complex case we need angle variables.  Will the 
result depend upon the choice of these angle variables?  
It appears that we can just take $\theta_i$ to be the angle for 
$x_i$.